\documentclass[ejs,preprint]{imsart}

\usepackage{amsfonts}
\usepackage{mathtools} 
\usepackage{amssymb, amsmath}
\usepackage{multicol}
\usepackage{multirow}
\usepackage[mathscr]{eucal}
\usepackage{graphics}
\usepackage{graphicx}
\usepackage{verbatim}
\usepackage{color}
\usepackage{float}

\usepackage{color}

    \newcommand{\argmin}{\mathop{\rm argmin}}

    \newcommand{\R}{{\mathbb R}}
    \newcommand{\E}{{\mathbb E}}

    \newcommand{\wh}{\widehat}
    \newcommand{\X}{\mathfrak{X}}

    \newtheorem{theorem}{Theorem}[section]
    \newtheorem{lemma}{Lemma}[section]

    \newtheorem{corollary}{Corollary}[section]

    \newtheorem{remark}{Remark}[section]
    \numberwithin{equation}{section}

\allowdisplaybreaks

\doi{10.1214/154957804100000000}
\pubyear{0000}
\volume{0}
\firstpage{0}
\lastpage{0}

\startlocaldefs
\endlocaldefs

\begin{document}

\begin{frontmatter}

\title{Nonparametric Confidence Regions for Level Sets: Statistical Properties and Geometry}
\runtitle{Confidence Regions for Level Sets}

\begin{aug}
  \author{Wanli Qiao\thanksref{t2}\ead[label=e1]{wqiao@gmu.edu}}

  \address{ Department of Statistics\\
George Mason University\\
4400 University Drive, MS 4A7\\
Fairfax, VA 22030\\
           \printead{e1}}

  \author{Wolfgang Polonik\thanksref{t3}
  \ead[label=e3]{wpolonik@ucdavis.edu}}

  \address{    Department of Statistics\\
    University of California\\
    One Shields Ave.\\
    Davis, CA 95616-8705\\
          \printead{e3}}

  \thankstext{t2}{Partially supported by NSF grant DMS 1821154.}
\thankstext{t3}{Partially supported by NSF grant DMS 1107206.}
  \runauthor{W. Qiao and W. Poloink}

\end{aug}

\begin{abstract}
This paper studies and critically discusses the construction of nonparametric confidence regions for density level sets. Methodologies based on both vertical variation and horizontal variation are considered. The investigations provide theoretical insight into the behavior of these confidence regions via large sample theory. We also discuss the geometric relationships underlying the construction of  horizontal and vertical methods, and how finite sample performance of these confidence regions is influenced by geometric or topological aspects. These discussions are supported by numerical studies. 
\end{abstract}




\tableofcontents

\end{frontmatter}

\section{Introduction}
For a density function $f$ on $\mathbb{R}^d,\; d\geq 1,$ we define the superlevel set of $f$ at level $c$ and the corresponding isosurface or contour as 
$$\mathcal{L} = \{x\in\mathbb{R}^d: \; f(x)\geq c\}\quad\text{and}\quad \mathcal{M} = \{x\in\mathbb{R}^d: \; f(x)= c\},$$ 
respectively. The dependence on the level $c$ is suppressed in our notation, for $c$ is fixed throughout the manuscript. We will study methods for constructing confidence regions for ${\cal M}$ and ${\cal L}$ based on an iid sample $X_1, X_2, \cdots, X_n$  from $f$. While new methods for constructing confidence regions are proposed below, and a comparison of these and other existing methods is provided, the overarching goal of this work is to provide a critical discussion of the advantages and disadvantages of the various existing methods. This will, among others, also lead to insight about the influence of geometry on the statistical performance of the confidence regions.\\  

The estimation of level sets (or isosurfaces) has received quite some interest in the literature. For some earlier work in the context of density level set estimation see Hartigan (1987), Polonik (1995), Cavalier (1997), Tsybakov (1997), Walther (1997). There are relations to density support estimation, (e.g. Cuevas et al. 2004, Cuevas 2009), clustering (e.g. Cuevas et al. 2000, Rinaldo et al. 2010), classification (e.g., Mammen and Tsybakov, 1999, Hall and Kang, 2005, Steinwart et al. 2005, Audibert and Tsybakov 2007), anomaly detection (e.g. Breuning et al. 2000, Hodge and Austin 2004), and more. Bandwidth selection for level set estimation is considered by Samworth and Wand (2010), and Qiao (2018a). Applications of level set estimation exist in many fields, such as astronomy (Jang 2006), medical imaging (Willett and Nowak, 2005), geoscience (Sommerfeld et al. 2015), and others. See Mason and Polonik (2009) and Mammen and Polonik (2013) for further literature. Level sets of functions also play a central role in topological data analysis, in particular in `persistent homology', where the topological properties of level sets as a function of the parameter $c$ are used for statistical analysis; e.g. see Fasy et al. (2014).\\

All the methods discussed in this paper are based on a kernel estimator 
\begin{align*}
\wh f(x)=\frac{1}{nh^d}\sum_{i=1}^nK\Big(\frac{X_i-x}{h}\Big), \;\; x\in\mathbb{R}^d,
\end{align*}
and its derivatives. Here $K$ is a $d$-dimensional kernel, and $h > 0$ is the bandwidth. As estimators for the superlevel sets or the corresponding contours, we are considering plug-in estimators given by
$$\widehat{\mathcal{L}}=\{x\in\mathbb{R}^d: \widehat f(x) \geq c\} \quad\text{and}\quad \widehat{\mathcal{M}}=\{x\in\mathbb{R}^d: \widehat f(x)=c\}.$$ 
It is well-known that $\widehat f(x)$ is a biased estimator. For a twice continuously differentiable kernel $K$, we will also consider a de-biased version of the kernel estimator, given by  $\wh f^{\,bc} = \widehat f - \widehat \beta$ with 
\begin{align*}
\wh \beta(x) = \frac{1}{2} h^2\int u_1^2K(u)du \sum_{j=1}^d  \frac{\partial^2}{\partial x_j \partial x_j} \wh f_l (x),
\end{align*}
where $\wh f_l$ denotes a kernel density estimator using bandwidth $l$, which can be different from $h$.  See Chen (2017) for a recent study of the de-biased estimator.\\

In some of the recent related literature (see, e.g., Chen et al., 2017), the bias is ignored entirely, and the target is redefined as a `smoothed' version of the superlevel set, or its contour, given by
$$\mathcal{L}^E = \{x\in\mathbb{R}^d: \; \E\wh f(x)\geq c\}\quad\text{and}\quad \mathcal{M}^E = \{x\in\mathbb{R}^d: \; \E\wh f(x)= c\},$$ 
respectively, where here and in what follows, we use the superscript $E$ to indicate that the definition is based on an `expected quantity'. We will also consider this target for some of our methods. \\[5pt]
For an interval $[a,b] \subset {\mathbb R},$ and a real valued function $g$, we denote by $g^{-1}[a,b] $ the pre-image of $[a,b]$ under $g$. For $a=b$, we also use the standard notation $g^{-1}(a)$ for the pre-image. Two different types of confidence regions are considered. One of them is based on {\em vertical variation}. Such confidence regions for ${\cal M}$ or ${\cal M}^E$ are of the form
\begin{align}\label{vert-conf}
\widehat{C} = \widehat{f}^{-1}\big[ c - \widehat{a}_n , c + \widehat{a}_n \big].
\end{align}
Corresponding lower and upper confidence regions for ${\cal L}$ (and ${\cal L}^E$) will be of the form $\widehat{f}^{-1}\big[c + \widehat{a}_n, \infty \big)$ and $\widehat{f}^{-1}\big[ c - \widehat{a}_n, \infty \big),$ respectively. The crucial question then is, how to choose the quantity $\widehat{a}_n$ in order to achieve a good performance. The other type of confidence region considered here is based on {\em horizontal variation}. Such confidence regions for ${\cal M}$ (or ${\cal M}^E$) are of the form
\begin{align}\label{horiz-conf}
\widehat{C} = \{y \in {\mathbb R}^d:\,\|y - x\| \le \widehat{b}_n(x), \, x \in \widehat{f}^{-1}(c)\},
\end{align}
where the random variable $\wh b_n(x)$ either does not depend on $x$, i.e. $\wh b_n(x) = \wh b_n$, or $\wh b_n(x) = \wh a_n/\|\nabla \wh f(x)\|$, for some random variable $\wh a_n$ not depending on $x$. Once horizontal variation based confidence regions $\wh C$ for ${\cal M}$ are constructed, we can obtain corresponding confidence regions for ${\cal L}$: Upper bounds (sets) are obtained by adding $\wh C$ to $\wh {\cal L}$, and lower bounds are constructed by subtraction. Confidence regions for ${\cal M}^E$ and ${\cal L}^E$ are constructed similarly.\\[4pt]
 For a simple  heuristic underlying the choice of $\wh b_n(x) = \wh a_n /\|\nabla \wh f(x)\|$, consider the one-dimensional case. 
%
For small $\wh a_n$, we have $|\wh a_n| \approx |\wh f(x+\wh b_n(x)) - \wh f(x)| = |\wh f(x+\wh b_n(x)) - c|$ for $x\in\wh f^{-1}(c)$, so that $\wh a_n$ reflects vertical variation, while $\wh b_n(x)$ represents horizontal variation.  For more details see Section~\ref{Horizontal}.\\[4pt]
%
%
Geometrically, confidence regions based on horizontal variation as in (\ref{horiz-conf}), but with a constant $\wh b_n,$ consist of tubes of constant width put around the estimated contour $\widehat{\cal M}.$ In other words, the maximum (horizontal) distance to ${\cal M},$ or ${\cal M}^E,$ is being controlled. Confidence regions based on vertical variation as in (\ref{vert-conf}), as well as the regions based on horizontal variation with non-constant $\wh b_n(x),$ have varying width, which contain information about the slope of the density. Different constructions of $\widehat{a}_n$ and $\widehat{b}_n$  will be discussed. One of them is based on extreme value theory for Gaussian fields indexed by manifolds, and the others on various bootstrap approximations. Some of the constructions of the horizontal methods also involve the estimation of integral curves. \\

The same type of confidence regions as discussed above, can also be constructed with $\wh f(x)$ replaced by the de-biased density estimator $\wh f^{bc}(x)$. We note that in the literature, different approaches to the removal of bias effect in confidence bands or regions using kernel-type estimators. One is based on undersmoothing of the original kernel density estimator. This approach is not considered in some detail here (but see Remark~\ref{discuss} below). Alternatively, one can use explicit bias corrections, or the smoothed bootstrap. Both of these approaches are being discussed in our work. For further relevant recent literature see Chen (2017) and Calonico et al. (2018a).\\
 


Bootstrap confidence regions for a density superlevel set and/or isosurface based on vertical variation can be found in Mammen and Polonik (2013) and Chen et al. (2017). The latter also constructed a confidence set based on horizontal variation, and such constructions can also be found in Chen (2017). The confidence sets proposed here are compared to these methods. In a different setting Jankowski and Stanberry (2014) consider confidence regions for the expected value of a random sets (or its boundary) based on repeated observations of the random set.  \\[5pt]
In summary, the contributions of the manuscript are to
\begin{itemize}
\item[(i)] derive asymptotically valid confidence regions based on vertical variation using extreme value distributions of kernel estimators indexed by certain manifolds; 
\item[(ii)] use bootstrap methods to construct confidence regions based on vertical variation in order to improve finite sample performance of the confidence regions in (i);
\item[(iii)] derive asymptotically valid bootstrap confidence regions based on horizontal variation;
\item[(iv)] discuss geometric connections between the different constructions of confidence regions;
\item[(v)] provide numerical studies to compare the finite sample performance of the various confidence regions developed here and in the literature;
\item[(vi)] critically discuss advantages and disadvantages of the two different types of constructing confidence regions  (horizontal and vertical variation).
\end{itemize}

The theoretical result underlying (i), see Theorem~\ref{LevelSetConf}, provides a closed form of asymptotically valid confidence regions for superlevel sets and isosurfaces. To the best of our knowledge, this is the first result of this type, and it provides important qualitative insight into the underlying problem. Due to the well-known slow convergence properties of extreme value distributions, the construction of bootstrap confidence regions appear of higher practical relevance.\\

In a recent paper, Qiao and Polonik (2018) derived the asymptotic distribution of extrema of certain rescaled  non-homogeneous Gaussian fields. This result provides an important tool for the construction of our confidence regions based on large sample distribution theory. The paper Qiao and Polonik (2016) on ridge (or filament) estimation is also a source of inspiration for the work presented here. \\

The paper is organized as follows. Section \ref{vertical} is considering vertical variation based methods. We first present a result on the behavior of the coverage probability of asymptotic confidence regions for isosurfaces and levels sets, and then we construct a bootstrap-based  confidence region. Section~\ref{Horizontal} is dedicated to horizontal methods. Section \ref{Horizontal1} presents the construction of two confidence regions of the form (\ref{horiz-conf}) using integral curves. Section \ref{Horizontal2} discusses a horizontal bootstrap based confidence region related to the Hausdorff-distance based method of Chen et al. (2017). In Section \ref{Discussion} we discuss how the finite sample behavior of our methods is influenced by some geometric or topological properties of estimated superlevel sets etc. Simulation results presented in Section \ref{simulations} compare the various methods. Section~\ref{Conclusion} presents some concluding discussions. The proofs of the technical results are presented in Section \ref{Proofs}.
\section{Confidence regions based on vertical variation}\label{vertical}

The following notation is used throughout the manuscript. For a sequence $\gamma > 0,$ we let 
\begin{align} \label{beta-notation}
    \beta^{(k)}_{n,\gamma} = \gamma^2 + \sqrt{\frac{\log n}{n \gamma^{d + 2k}\,}},\quad \text{and}\quad \beta_{n,\gamma}^{(k),E} = \sqrt{\frac{\log n}{n \gamma^{d + 2k}\,}},\quad k = 0,1,2,3.
\end{align}
For a kernel density estimator based on a bandwidth $h,$ the quantity $\beta_{n,h}^{(0)}$ equals the rate of the uniform deviation from the density under standard assumptions (satisfied in our setting). In other words, we have uniform consistency of the kernel estimator if $\beta_{n,h}^{(0)} \to 0$ as $n \to \infty.$ The quantities $\beta_{n,h}^{(k)}, k = 1,2,3$ have the same interpretation when considering the kernel estimator of the $k$-th derivatives of the density. Similarly, $\beta_{n,h}^{(k),E}$ is the standard uniform rate of convergence of the kernel density estimator of the $k$-th derivative when centered at its expectation, $k = 0,\ldots,3$. 

\subsection{Confidence regions based on asymptotic distribution}
Our first result provides asymptotically valid confidence regions for the isosurface $\mathcal{M}$ and the superlevel set $\mathcal{L}$. Before formulating the theorem we introduce the underlying assumptions and some more notation.\\[5pt]
{\sc Assumptions:\\[5pt]}
%
%
(\textbf{F1}) The probability density $f$ has bounded, continuous derivatives up to fourth order. There exists $m>0$ such that $\int \|x\|^m f(x)dx<\infty$.\\[5pt]
(\textbf{F2}) There exist $\delta_0>0$ and $\epsilon_0>0$ such that $\|\nabla f(x)\| > \epsilon_0$ for $x\in \{x: c - \delta_0 \leq f(x) \leq c + \delta_0\}.$\\[5pt]
%
%
(\textbf{K}) The kernel function $K: \; \mathbb{R}^d \rightarrow \mathbb{R}$ is symmetric about zero, with support contained in $[-1,1]^d$, and is continuously differentiable up to fourth order.\\[5pt]
(\textbf{A}) Both $f$ and $K$ are continuously differentiable up to $d$-th order.\\[5pt]
%
(\textbf{H1})$^k$. $h\rightarrow0$, and $\beta_{n,h}^{(k),E}\rightarrow0$ as $n\rightarrow\infty$, where $k = 0$ or $2$. \\[5pt]
({\bf H2})$^k$ The bandwidth $l$ used for bias correction, satisfies $l \rightarrow 0$ and $\beta_{n,l}^{(2k),E}\rightarrow0$, where $k = 0$ or $2$. In addition, $(h/l)^{\frac{d+4}{2}}\log{n} \rightarrow 0$ and $\sqrt{nh^{d+4}\log{n}\,}\,l^2 \rightarrow 0$. 

\begin{remark} \label{remark}$\;$\\[3pt]
{\em
a) Assumption ({\bf F1}) can be weakened when only considering confidence regions for the smoothed isosurface ${\cal M}^E$. Only continuous second order derivatives of $f$ are needed in this case. For simplicity, we just use the stronger assumption ({\bf F1}) throughout the manuscript.\\[3pt]
b) Assumptions ({\bf F1}) (even the weakened version discussed in 1.) and ({\bf F2}) together imply that the $(d-1)$-dimensional isosurface $\mathcal{M}$ has positive reach (e.g. see Lemma 1 of Chen et al., 2017), which is a necessary condition for applying the result in Qiao and Polonik (2018). Assumptions on the reach (introduced by Federer 1959) are used in a number of studies of geometric properties of manifolds, etc. Also note that assumption ({\bf F2}) implies that the level $c>0$.\\[3pt]
%
%
c) It seems likely that the assumption of a non-negative kernel can be replaced by a higher-order kernel. This is because our focus is on the superlevel sets (or isosurfaces) with $c>0$ and regions there the density estimator is negative is irrelevant to our estimation. Even for the bootstrap method, technically we can bootstrap from a truncated and normalized density estimator. However, the truncation might incur further technical considerations.\\[3pt]
%
%
%
d) The smoothness requirements for $f$ and $K$ in assumption ({\bf A}) are needed to enable to use of the Rosenblatt transform (see Rosenblatt, 1976) in the proof of Theorem 2.1. It will not be needed for bootstrap based methods. This is why we did not combine assumptions {\bf (A)}, ({\bf F1}) and ({\bf K}).\\[3pt]
e) Assumptions (\textbf{H1})$^0$ and (\textbf{H2})$^0$, both very mild, are used in large sample confidence regions based on the vertical method. They guarantee the uniform consistency of the de-biased estimator $\wh f^{bc}=\wh f -\wh\beta$. Specifically, (\textbf{H1})$^0$ is used for $\wh f$ centered at its expectation, while (\textbf{H2})$^0$ is used for the bias correction part $\wh\beta$. We use stronger assumptions (\textbf{H1})$^2$ and (\textbf{H2})$^2$ for confidence regions based on the horizontal method. In particular, these assumptions guarantee the uniform consistency of the second derivatives of $\wh f^{bc}$.\\[3pt]
f) Under our assumptions, in the case of $d=1$, $\cal  M$ is a union of finitely many points, say $N$. Denote $\mathcal{M} = \{x_i, i=1,\cdots, N\}$.
}
\end{remark}

We need to introduce further notation. Let $\beta(x) = \mathbb{E} \wh{f} (x) - f(x)$ be the bias. Set $\|K\|_2^2 = \int K^2(u)du,$ and let $\mathscr{V}_{d-1}$ denote $(d-1)$-dimensional Hausdorff measure. Let $\wh{\mathscr{V}_{d-1}(\mathcal M)}$ be an estimator of $\mathscr{V}_{d-1}(\mathcal M)$. Some specific estimators will be given in Remark~\ref{discuss}c). For $d \ge 2$ and $s^2_K = \frac{\int \big[\frac{\partial}{\partial u_1} K(u)\big]^2du}{2\int K^2(u)du},$ we set 
\begin{align}\label{hatb-def}
\wh b(\alpha) &= \sqrt{2(d-1)\log{h^{-1}}}+\frac{1}{\sqrt{2(d-1)\log{h^{-1}}}}\bigg[z(\alpha)+\left(\frac{d}{2}-1\right)\log{\log{h^{-1}}}\nonumber\\
&\hspace{4cm}+\log\bigg\{\frac{(2d-2)^{d/2-1} s_K^{d-1}}{\sqrt{2}\pi^{d/2}}\wh{\mathscr{V}_{d-1}(\mathcal M)}\bigg\}\bigg].
\end{align}
where $0 < \alpha < 1$, and $z(\alpha) =  -\log[-\frac{1}{2} \log(1-\alpha)]$. The quantity $\wh b(\alpha)$ is based on the extreme value behavior of a Gaussian field indexed by ${\cal M}$ (see (\ref{ThetaExp}) given in Theorem~\ref{ProbabilityResult}). We further need the following, which, for $d \ge 2$, simply is a scaled version of $\wh b(\alpha):$ 
\begin{align}\label{def-hata}
\wh a^{(d)}_{1-\alpha} = 
\begin{cases}
\frac{ \wh b(\alpha)\sqrt{\|K\|_2^2c}}{\sqrt{nh^d}} & \text{for }\;d \ge 2\\
\frac{ \Phi^{-1}\big((1-\alpha)^{1/\wh N}\big)\sqrt{\|K\|_2^2c}}{\sqrt{nh^d}}& \text{for }\;d = 1,
\end{cases}
\end{align} 
where $\Phi$ is the standard normal cdf, and, for $d=1$, $\wh N$ is the cardinality of $\wh{\cal M}$ (cf. Remark~\ref{remark}f)\,). Note that when $d\geq2$, $\wh a^{(d)}_{1-\alpha}$ has a typical structure appearing in confidence bands for probability density or regression functions. See, for example, the main result of Bickel and Rosenblatt (1973). When $d=1$, $\wh a^{(d)}_{1-\alpha}$ corresponds to a quantile value of the maximum of a Gaussian mixture model because $\cal M$ is a collection of separated points under our assumptions.  
%
%
With this notation, our first confidence interval based on vertical variation, and using the de-biased estimator of the underlying density, is defined as 
$$\wh C_{n,1}(1-\alpha) = (\wh f^{\,bc})^{-1}\left[c -  \wh a^{(d)}_{1-\alpha},c +  \wh a^{(d)}_{1- \alpha}\right].$$
%
%
\begin{theorem}\label{LevelSetConf}
Suppose that (\textbf{F1}), (\textbf{F2}), (\textbf{K}), ({\textbf A}), ({\textbf H1})$^0$ and ({\textbf H2})$^0$ hold. Let $0 < \alpha < 1.$ If $\wh{\mathscr{V}_{d-1}(\mathcal M)}$ is a consistent estimator for $\mathscr{V}_{d-1}(\mathcal M)$, then we have, 
\begin{align}
\lim_{n\rightarrow \infty}\mathbb{P}\Big( \mathcal{M} \subset \wh C_{n,1}(1-\alpha)  \Big) =1-\alpha.\label{ContourConf}
\end{align}

%
\end{theorem}

%
%
%
%
%

An asymptotic confidence region for the superlevel set $\mathcal{L}$ is given by the following upper and lower bounds:
\begin{align*}
\wh C_{n,1}^{-}(1-\alpha) = 
(\wh f^{\,bc})^{-1}\Big[c -  \wh a^{(d)}_{1- \alpha} , +\infty \Big)
\end{align*}
and
\begin{align*}
\wh C_{n,1}^{+}(1-\alpha) = 
(\wh f^{\,bc})^{-1}\Big[c +  \wh a^{(d)}_{1- \alpha} , +\infty \Big).
\end{align*}

\begin{corollary}\label{UpperLSConf}
Suppose that (\textbf{F1}), (\textbf{F2}),  (\textbf{K}), ({\textbf A}), ({\textbf H1})$^0$ and ({\textbf H2})$^0$ hold. Also suppose $\beta_{n,h}^{(1),E}\rightarrow0$ and $\beta_{n,l}^{(3),E}\rightarrow0$ as $n\rightarrow\infty$. Let $0 < \alpha < 1.$ If $\wh{\mathscr{V}_{d-1}(\mathcal M)}$ converges to $\mathscr{V}_{d-1}(\mathcal M)$ in probability, then we have
\begin{align*}
\lim_{n\rightarrow \infty}\mathbb{P}\Big( \wh C_{n,1}^{+}(1-\alpha) \subset \mathcal{L} \subset \wh C_{n,1}^{-}(1-\alpha)  \Big) =1-\alpha.
\end{align*}
\end{corollary}
\begin{remark}\label{discuss}  $\;$
{\em 
a) Notice that the construction of $\wh C_{n,1}$ involves the choice of two bandwidths, $h$ and $l$. This is the case for all the confidence regions considered in this paper that are based on $\wh f^{\,bc}$.\\[5pt]
%
b) Constructing $\wh C_{n,1}(1-\alpha)$ with $\wh f$ rather than $\wh f^{\,bc}$ also results in an asymptotically valid confidence set for ${\cal M},$ provided undersmoothing is being used to handle the bias. In this case, a sufficient condition for consistency of the coverage probability is that the bandwidth satisfies $\beta_{n,h}^{(2),E} \rightarrow 0$. When this assumption holds, the stochastic term $\sup_{x\in\mathcal{M}} |\wh f(x) - \mathbb{E}\wh f(x) |$ asymptotically dominates the bias term $\sup_{x\in\mathcal{M}} |\mathbb{E}\wh f(x) - f(x)|$, so that the latter can be ignored in the proof (see Hall, 1993). \\[5pt]
%
c) We discuss two choices for the estimator $\wh{\mathscr{V}_{d-1}(\mathcal M)}$. Let $\lambda$ be the $d$-dimensional Lebesgue meaure, and let $P_n = n^{-1}\sum_{i=1}^n \delta_{X_i}$ denote the empirical probability measure, where $\delta_{x}$ denotes Dirac measure in $x$. Let $\mathscr{A}$ be the class of compact sets with a positive reach bounded away from zero such that $\mathcal{L}\in\mathscr{A}$. One option for the consistent estimator of $\mathscr{V}_{d-1}(\mathcal M)$ is given by $\mathscr{V}_{d-1}(\partial \wh{\mathcal L})$, where
\begin{align*}
\wh{\mathcal L} = \argmin_{A\in\mathscr{A}} [P_n(A) - c\lambda(A)].
\end{align*}
The consistency of this estimator is shown in Proposition 3 in Cuevas et al. (2012). There, consistency is derived in terms of outer Minkowski content, which, under our assumptions, is equivalent to the consistency using Hausdorff measure (see Corollary 1 in Ambrosio et al., 2008). Efficiently computing $\wh{\mathcal L}$ is challenging. 

Another estimator for $\mathscr{V}_{d-1}(\mathcal M)$ is given by  $\mathscr{V}_{d-1}(\wh{\mathcal M})$ with $\wh{\mathcal M} = \wh{f}^{-1}(c).$ The convergence rate and asymptotic normality of this estimator in the context of surface integral estimation are shown in Theorem 1 in Qiao (2018b), where additional assumptions are imposed, in particular on the speed of convergence of $h$.  \\[5pt]
d) The $(d-1)$-dimensional Hausdorff measure of $\mathcal{M}$ and its estimators come into play because (i) the distribution of the statistic related to the above confidence regions can be approximated by that of the extreme value of certain Gaussian random fields indexed by the level set; and (ii) the latter, in turn, is related to an integral over $\mathcal{M}$ with respect to the $(d-1)$-dimensional Hausdorff measure.  In our case, the integrand of this integral turns out to be a constant - see Theorem 6.1, and thus we obtain the volume of the isosurface. In fact, it is well known that the probability of the extreme value of a locally stationary Gaussian random fields exceeding a large level is asymptotically proportional to the volume of the index set (locally speaking). See, for example, Chapter 2 of Piterbarg (1996). Surface integrals have appeared in the context of level-set estimation before in Cadre (2006). There, however, a first order asymptotics (consistency) is considered, using the set-theoretic measure of symmetric difference $d({\cal L}, \wh {\cal L})$. On a very heuristic level, the fact that a surface integral comes in here can be understood by the fact that the variance of the fluctuations of $d({\cal L}, \wh {\cal L})$ is of the order $a_n = 1/nh^d,$ which is inherited from the fluctuations of the density estimator, and by then approximating $d({\cal L}, \wh {\cal L})$ by a constant times ${\rm vol}\big(f^{-1}[c+a_n, c - a_n]\big),$ the Lebesgue measure of $f^{-1}[c+a_n, c - a_n].$ Lebesgue's theorem gives that $a^{-1}_n{\rm vol}\big(f^{-1}[c+a_n, c - a_n]\big)$ converges to a surface integral over the boundary ${\cal M}$. In other words, the technical reason for this surface integral to appear in Cadre (2006) is different from why it appears in our context.
}
\end{remark}

\subsection{Bootstrap confidence regions}\label{BootstrapVertical}
Bootstrap confidence regions based on vertical variation of the kernel density estimate have been constructed in Mammen and Polonik (2013) and in Chen et al. (2017). They are based on a bootstrap approximation of quantiles of statistics of the form
$$T(D_n) = \sup_{x \in D_n}\big|\sqrt{nh^d\,}\, (\wh f(x) - f(x))\big|,\text{ \rm or }T^E(D_n) = \sup_{x \in D_n}\big|\sqrt{nh^d\,}\, (\wh f(x) - {\mathbb E}\wh f(x))\big|,$$
where $D_n$ is such that (asymptotically) $ {\cal M} = f^{-1}(c) \subset D_n.$ The confidence regions considered in the literature differ in the choice of the set $D_n.$ While Mammen and Polonik (2013) propose to use $D_n = \{c - \epsilon_n \le f(x) \le c + \epsilon_n\}$ for some appropriate choice of $\epsilon_n$ tending to zero, as $n \to \infty$, Chen et al. simply use $D_n = {\mathbb R}^d.$ Here we are using the smallest such set $D_n = {\cal M}.$  Clearly, the statistics are stochastically ordered in terms of the size of the set $D_n$. Thus, the choice $D_n =  {\cal M} $ leads to the confidence set that is the smallest among the three. Of course the coverage still needs to be investigated. However, if the corresponding bootstrap approximations of the distributions of $T(D_n)$ (and $T^E(D_n)$) work similarly well, then using the statistic $T({\cal M})$ or $T^E({\cal M})$, respectively, can be expected to be a good choice for the construction the bootstrap based confidence sets based on the vertical variation.\\[5pt]
Our construction is as follows. Let $X_1^*,\ldots,X_n^*$ be a sample drawn from a kernel density estimator $\wh f_g$ using $X_1,\ldots,X_n,$ and bandwidth $g>0$. Let $\wh f^*(x)$ be the kernel density estimate using $X_1^*,\ldots,X_n^*$ and bandwidth $h$. Let $\wh f^{*,E}(x) = \mathbb{E}^* \wh f^*(x)$, where we use $\mathbb{E}^*$ to indicate the expectation with respect to $\wh f_g$. For $0 < \alpha < 1$, let $\widehat c_{1-\alpha}^{*,E}$ be the $(1-\alpha)$-quantile of the distribution of $\sup_{x\in\wh{\mathcal{M}}}|\widehat f^*(x) - \wh f^{*,E}(x)|$, and let $\widehat c_{1-\alpha}^*$ be the corresponding quantile of $\sup_{x\in\wh{\mathcal{M}}}|\widehat f^*(x) - \widehat f_g(x)|$. We now define our bootstrap confidence regions for ${\cal M}$ and ${\cal M}^E$, respectively, as
\begin{align}\label{Cn2}
&\wh C_{n,2}^*(1-\alpha) =\wh f^{-1}[c - \widehat c_{1-\alpha}^*, c + \widehat c_{1-\alpha}^*] \nonumber\\
\text{ \rm and } &\wh C_{n,2}^{*,E}(1-\alpha) = \wh f^{-1}[c - \widehat c_{1-\alpha}^{*,E}, c + \widehat c_{1-\alpha}^{*,E}].
\end{align}
We also define the following sets to construct bootstrap confidence regions for ${\cal L}$ and ${\cal L}^E$, respectively.
\begin{align*}
&\wh C_{n,2}^{*,-} (1-\alpha)=\wh f^{-1}[c - \widehat c_{1-\alpha}^*, +\infty), \quad \wh C_{n,2}^{*,+} (1-\alpha) =\wh f^{-1}[c + \widehat c_{1-\alpha}^*,+\infty), \\
 &\wh C_{n,2}^{*,E,-} (1-\alpha)= \wh f^{-1}[c - \widehat c_{1-\alpha}^{*,E}, \infty), \quad\text{\rm and}\quad C_{n,2}^{*,E,+} (1-\alpha)= \wh f^{-1}[c + \widehat c_{1-\alpha}^{*,E}, +\infty).
\end{align*}
Below we show that this (and other) confidence region is asymptotically exact, and we derive rates of convergence for the coverage probability. These rates of convergence have a somewhat complex appearance, which we first explain from a high level perspective.
\vspace*{0.2cm}

{\bf Structure of the rates of convergence for the coverage probabilities of bootstrap based confidence sets:} The derivation of the following somewhat complex looking rates of convergence of the coverage probabilities are all based on Lemma~\ref{MP}, which is a slight reformulation of a result of Mammen and Polonik (2013). Based on this result, the rates are of the form $O\big(\sqrt{\delta_n} + \tau_n),$ where $\delta_n$ and $\tau_n$ are derived as follows: Let $Z_n = \sup_{x \in A}\big|\wh g(x) - \overline{g}(x) \big|,$ and let $Z_n^*$ be a bootstrap version, where $\wh g(x)$ is one of the density estimators considered, and $\overline{g}(x)$ is some centering; the set $A$ in the supremum is either $\cal M$ or ${\cal M}^E.$ Then, we derive $\delta_n$ and $\tau_n$ by showing that for some sequence of positive real numbers $\gamma_n$, we have $P\big(\big|Z_n - Z_n^*\big| > \gamma_n\big) \le \delta_n$ and $\sup_{t \in \R}P\big(Z_n \in [t, t+\gamma_n)\big) \le \tau_n.$ Both of $\delta_n$ and $\tau_n$ are themselves comprised of a sum of various terms. In fact, in our applications of this result, $\gamma_n$ is chosen such that $\sqrt{\delta_n}$ is negligible, and we have $\tau_n = \Psi_n(\gamma_n)$, with 
\begin{align}\label{Neumann}
\Psi_n(\gamma) = \gamma \sqrt{nh^d\log{n}} + h\log{n} + \log{n}\left( \sqrt{\beta_{n,h}^{(0),E}} + \sqrt{h^d\log{n}}\right),\quad \gamma > 0.
\end{align}
That, in our case, $\tau_n$ has this particular form follows from a result by Neumann (1998). The rates $\gamma_n$ that make $\sqrt{\delta_n}$ negligible will follow, respectively, from strong approximation results in Neumann (1998) and some modification of Neumann's result given in Mammen and Polonik (2013).\\[3pt]
Note further that in each of the construction approaches considered below, we obtain the same rates of convergence of the coverage probabilities for both confidence regions for ${\cal M}$ and confidence regions for ${\cal L}$ (and similarly for ${\cal M}^E$ and ${\cal L}^E$). The reason for this is as follows. Let $p_{n,\cal M}$ and $p_{n,\cal L}$ denote these coverage probabilities based on one approach (e.g. the left-hand sides of (\ref{confm2}) and (\ref{confl2}), respectively. By construction, we have $p_{n,\cal L} \le p_{n,\cal M}$. More precisely, it is shown in the proof of Corollary~\ref{UpperLSConf} that $p_{n,\cal L} = p_{n,\cal M} \cdot q_n$, where $0 \le 1 - q_n = O(\delta_n)$ with $\delta_n = o(1).$ Now, if $p_{n,\cal M} = 1-\alpha + O(\tau_{n}),$ then we obtain $p_{n,\cal L} = (1-\alpha + O(\tau_{n})) (1 - O(\delta_n)) = (1-\alpha) + O\big(\max(\tau_n,\delta_n)\big)$, and under our respective assumptions, $O\big(\max(\tau_n,\delta_n)\big) = O(\tau_n),$ showing that both $p_{n,\cal M}$ and $p_{n,\cal L}$ converge to $(1-\alpha)$ at the same speed. 

%
%
%
\begin{theorem}\label{smooth-bootstrap-method}
Suppose that (\textbf{F1}), (\textbf{F2}) and (\textbf{K}) hold. Let $0 < \alpha < 1,$ and
%
%
\begin{align}
\label{gammanE}
\gamma_n^E &= \Big(\sqrt{\beta_{n,g}^{(0)}\,} +  \beta_{n,h}^{(1),E }\,\Big) \beta_{n,h}^{(0),E},\\
\gamma_n &= \gamma_n^E + \Big( \beta_{n,g}^{(2)}  + \beta_{n,h}^{(1)} \Big) h^2.
\label{gamman}
\end{align}
(a) If $\Psi_n(\gamma_n^E)=o(1),$ then we have
\begin{align}
&\mathbb{P} \left( \mathcal{M}^E \subset \wh C_{n,2}^{*,E} (1-\alpha)\right) = (1-\alpha) + O\big(\Psi_n(\gamma_n^E)\big),\quad \text{and} \label{confme2}\\
&\mathbb{P} \left( \wh C_{n,2}^{*,E,+} (1-\alpha) \subset \mathcal{L}^E \subset \wh C_{n,2}^{*,E,-} (1-\alpha)\right) = (1-\alpha) + O\big(\Psi_n(\gamma_n^E)\big).\label{confle2}
\end{align}
(b) If $\Psi_n(\gamma_n)=o(1)$ as $n\rightarrow\infty,$ then we have
%
\begin{align}
&\mathbb{P} \left( \mathcal{M} \subset \wh C_{n,2}^* (1-\alpha)\right) = (1-\alpha) + O\big(\Psi_n(\gamma_n)\big),\quad \text{and}\label{confm2}\\
&\mathbb{P} \left( \wh C_{n,2}^{*,+} (1-\alpha) \subset \mathcal{L} \subset \wh C_{n,2}^{*,-} (1-\alpha)\right) = (1-\alpha) + O\big(\Psi_n(\gamma_n)\big).\label{confl2}
\end{align}
\end{theorem}

\begin{remark} \label{vert-boot}

{\em a) The set $\wh C_{n,2}^{*,E}$ can be also used as a confidence region for $\mathcal{M}$ (not just for $\mathcal{M}^E$) if the bandwidth is chosen to be of smaller order than the optimal bandwidth, to make the bias negligible (undersmoothing). In practice, the choice of undersmoothing bandwidth might be difficult to determine. We do not pursue the theoretical justification for $\wh C_{n,2}^{*,E}$ as a confidence region for $\mathcal{M}$, the numerical performance of which, however, is shown in the simulation section.\\[5pt]
b) The quantity $\beta_{n,g}^{(2)}$ appears in $\gamma_n$, because the second partial derivatives appear in the bias of density estimation and need to be estimated in our proof. Note that we are not using the de-biased density estimator in this theorem.\\[5pt]
c)  For $d\geq 2$, if we choose both $g$ and $h$ to be of the standard optimal rates, i.e. $g=h= {\rm const.}\,n^{-1/(d+4)}$, then
\begin{align*}
\rho_n^E = O(h(\log{n})^{3/2}) = O(n^{-1/(d+4)}(\log{n})^{3/2} ).
\end{align*}
The rate in Chen et al. (2017) is given by $(nh^d)^{-1/8}(\log{n})^{7/8}$. When $h$ is chosen as the standard optimal bandwidth for density estimation, i.e. $h$ is of the exact order $n^{-1/(d+4)}$, the rate given in Chen et al. (2017) becomes $n^{-1/(2d+8)}(\log{n})^{7/8}$, which is slower than $\rho_n^E$. Note that even if the rate in Chen et al. (2017) is $(nh^d)^{-1/6}(\log{n})^{7/6}$ as claimed in this paper, it is still slower than $\rho_n^E$ when using the optimal bandwidth.\\

For $d\geq 2$, if we again choose $h=O(n^{-1/(d+4)})$, then
\begin{align*}
\rho_n = O\left( n^{-1/(d+4)}(\log{n})^{5/4} + \beta_{n,g}^{(2)}\sqrt{\log{n}} + \sqrt{\beta_{n,g}^{(1)}} \log{n} \right).
\end{align*}
With $g = O(n^{-1/(d+6)})$, $\rho_n$ is of the order of $n^{-1/(d+6)} \log{n}$. 
}
\end{remark}

\subsubsection{A bootstrap confidence region based on explicit bias correction}

Introducing a new bandwidth $g$ and bootstrapping from $\hat f_g,$ as we did above, can be viewed a method of bias correction (see page 208, Hall 1992). This allows us to construct a confidence region for ${\cal M}$ (rather than just for ${\cal M}^E$). We can also construct a confidence region for $\mathcal{M}$ using an explicit bias correction by
\begin{align*}
\wh C_{n,3}^*(1-\alpha) = (\wh f^{\,bc})^{-1} [c - \widehat c_{1-\alpha}^{*,E}, c + \widehat c_{1-\alpha}^{*,E}].
\end{align*}
For confidence regions for ${\cal L}$, we define
\begin{align*}
&\wh C_{n,3}^{*,-}(1-\alpha) = (\wh f^{\,bc})^{-1} [c - \widehat c_{1-\alpha}^{*,E}, +\infty),\quad \text{ and } \\
& \wh C_{n,3}^{*,+}(1-\alpha) = (\wh f^{\,bc})^{-1} [c + \widehat c_{1-\alpha}^{*,E}, +\infty).
\end{align*}

We have the following result:

\begin{theorem}\label{conf-region-explicit-bias}
Suppose that (\textbf{F1}), (\textbf{F2}), and (\textbf{K}) hold. Let $0 < \alpha < 1.$ Let 
\begin{align}\label{gammanbc} 
\gamma_n^{bc}=  \Big( \sqrt{\beta_{n,g}^{(0)}\,} + \beta_{n,h}^{(1)}\Big)\beta_{n,h}^{(0)} + h^2\beta_{n,l}^{(2)}.
\end{align}
If $\beta_{n,l}^{(2)}=o(1)$ and $\Psi_n\big(\gamma_n^{bc}\big)=o(1)$ as $n\rightarrow\infty$, then we have 
\begin{align}\label{cn3confm}
\mathbb{P} \left(  \mathcal{M} \subset \wh C_{n,3}^*(1-\alpha) \right) = (1-\alpha) + O\big(\Psi_n(\gamma_n^{bc})\big).
\end{align}
If we further assume $\beta_{n,l}^{(3)} =o(1)$ as $n\rightarrow\infty$, then
\begin{align}\label{cn3confl}
\mathbb{P} \left( \wh C_{n,3}^{*,+} (1-\alpha) \subset \mathcal{L} \subset \wh C_{n,3}^{*,-} (1-\alpha)\right) = (1-\alpha) + O\big(\Psi_n(\gamma_n^{bc})\big).
\end{align}
\end{theorem}

\begin{remark}{\em 
The constructions of both $\wh C_{n,2}^{*,E}$ and $\wh C_{n,3}^*$ are using a quantile, $\wh c^{*,E}_{1-\alpha},$ that is ignoring the bias. Nevertheless, $\wh C_{n,2}^{*,E}$ gives a confidence region for the smoothed isosurface ${\cal M}^E$, while $\wh C_{n,3}^*$ is a confidence region for ${\cal M}$. This is so, because one of them, $\wh C_{n,3}^*,$ is based on the de-biased estimator, while $\wh C_{n,2}^{*,E}$ is not.
Heuristically, this can be understood by writing $\wh f^{\,bc}(x) - f(x) = \wh f(x) - \mathbb{E}\wh f(x) + \beta(x) - \wh\beta(x).$ One can see that the bias correction in the density will adjust for the bias that is present in the quantile $\wh c^{*,E}_{1-\alpha}.$ \\

For $d\geq 2$, if we choose the optimal bandwidth $g=h=O(n^{-1/(d+4)})$, and $l=O(n^{-1/(d+8)})$ (which is the order of the optimal bandwidth for estimating the second derivatives), then
\begin{align*}
\rho_n^{bc} = O(h(\log{n})^{3/2}) = O(n^{-1/(d+4)}(\log{n})^{3/2} ).
\end{align*}
Compared to the rate $\rho_n$ given in Remark~\ref{vert-boot} this rate is faster. However, the fact that the construction of $\wh C_{n,3}^*$ involves the choice of three bandwidths might be a challenge in practice.
}
\end{remark}
\section{Confidence regions based on horizontal variation}\label{Horizontal}

Various confidence regions for ${\cal M}, {\cal M}^E, {\cal L}$ and ${\cal L}^E$ based on horizontal variation will be derived in this section. The geometric link between horizontal and vertical variation based confidence regions is, for obvious reasons, given by the gradient. Let $x \in {\cal M}$ and $\wh x \in \wh {\cal M}$, such that $\wh x$ is close to $x$, then we obviously have
\begin{align}\label{vert-horiz}
\wh f(x) - f(x) = \wh f(x) - \wh f(\wh x) \approx \nabla\wh f(\wh x)^T \big(x - \wh x\big),
\end{align}
and when $\wh x$ is chosen such that $x-\wh x$ is approximately perpendicular to $\wh{\cal M}$,
%
\begin{align}\label{vert-horiz2}
 \| x - \wh x \| \approx \frac{|\wh f(x) - f(x)|}{\|\nabla\wh f(\wh x)\|}.
\end{align}
Different ways of choosing $\wh x$ will give rise to different types of horizontal variation based confidence regions. For example, $\wh x$ can be a projection point of $x$ onto $\wh{\cal M}$ or it can be chosen such that there exists a gradient integral curve connecting $x$ and $\wh x$.\\[3pt]
While the above confidence regions based on vertical variation are using approximations to quantiles of the distribution of $\sup_{x \in {\cal M}}|\wh f(x) - f(x)|$, the horizontal variation based confidence regions will use estimated quantiles of the quantity $\sup_{x \in {\cal M}}\big|\nabla\wh f(\wh x)^T \big(x - \wh x\big)\big|$. The latter methods are not purely horizontal variation based, as they involve the adjustment of $\wh x - x$ by the gradient of $f$. Nevertheless, since they are explicitly using the differences $x - \wh x$, we still call them methods based on horizontal variation. \\[3pt]
Confidence regions based on horizontal variation (without estimating the gradient), have been constructed in Chen et al. (2017) based on the Hausdorff distance. The approaches considered in our work are asymptotically equivalent to the method proposed in Chen et al. (2017), but in practice the confidence sets are different. The different constructions also provide additional insight into the underlying geometry.\\[5pt]
In the following we introduce novel horizontal variation based methods, based on estimating quantiles via both the asymptotic distribution and the bootstrap. Instead of using standard bootstrap as in Chen et al. (2017), we adopt smoothed bootstrap. 
We would like to note that rather than simply providing alternative methods for confidence regions, the constructions and the discussion of the performance of the resulting confidence regions are meant to provide insight into the effects of geometric aspects of the underlying probability density on the performance of confidence regions.

\subsection{Methods based on integral curves}\label{Horizontal1} The approach discussed here is based on the construction of $\widehat{x}$ (cf. (\ref{vert-horiz})) using integral curves driven by the (scaled) gradient field and their relation to level set. This relation will be discussed first.\\[5pt]
{\sc Integral curves and level sets:} For any $x_0\in\mathbb{R}^d$, let $\{\X_{x_0}(t), t\in\mathbb{R}\}$ be the integral curve driven by the scaled gradient of $f,$ starting from $x_0$, defined by the equation, 
\begin{align*}
\frac{d\X_{x_0}(t)}{dt}=\frac{\nabla f (\X_{x_0}(t))}{\|\nabla f (\X_{x_0}(t))\|^2},\;\; \X_{x_0}(0)=x_0,
\end{align*}
where we assume that $\|\nabla f (\X_{x_0}(t))\| \ne 0$ (cf. Assumption ({\bf F2})). For $d=1$, $\nabla f$ is understood to mean $f^\prime$. The reason for choosing the scaled gradient $\frac{\nabla f}{\|\nabla f\|^2}$ as a driving vector field, rather than $\nabla f$ itself, is based on the following convenient property. Suppose that $t > 0$, and set $I_0(t) = [0,t]$. Then, by the fundamental theorem of calculus for line integrals, we have for any $x_0\in \cal M$ and $t\in\mathbb{R}$ such that $\|\nabla f (x)\| \ne 0$ for all $x \in \{s: c \le f(s) \le c+t\},$ 
\begin{align}
f (\X_{x_0}(t)) - f (x_0) &= \int_{\{\X_{x_0}(s): s\in I_0(t)\}}\nabla f (\textbf{r}) \cdot d\textbf{r}\label{eq5}\\
&=\int_{I_0(t)}\nabla f ( \X_{x_0}(s)) \cdot \frac{\nabla f( \X_{x_0}(s))}{\|\nabla f( \X_{x_0}(s))\|^2} ds\nonumber\\
&=\int_{I_0(t)} ds = t,\label{gradientpathproperty}
\end{align}
where the right hand side of (\ref{eq5}) is a line integral over the trajectory $\{\X_{x_0}(s): s\in I_0(t)\}$, and ``$\cdot$'' represents dot product between vectors. 
By adopting the current scaled vector field, we can make sure the height increase (or decrease) in the density is exactly the amount of ``time'' needed to travel. In other words, if two particles start from any two points $x_1, x_2 \in \mathcal{M}$, then, after following the integral curves $\X_{x_1}(\cdot)$ and $\X_{x_2}(\cdot)$ respectively for time $t$ (which can be negative), both of these two particles will arrive at $\mathcal{M}_{c+t}$. The same holds for $t \le 0$, by using the convention to start integration at $0$.\\[5pt] 
%
%
Let the plug-in estimators of $\nabla f(x)$ and $\X_{x_0}(t)$ based on the kernel density estimate be denoted by $\nabla \widehat f(x)$ and $\widehat \X_{x_0}(t)$, respectively, where the latter is the solution of the differential equation
\begin{align}\label{est-int-curve}
\frac{d\widehat\X_{x_0}(t)}{dt}=\frac{\nabla \widehat f (\widehat\X_{x_0}(t))}{\|\nabla \widehat f (\widehat\X_{x_0}(t))\|^2},\;\; \widehat\X_{x_0}(0)=x_0.
\end{align}
%
%
%
We use the notation $\wh\X^{bc}_{x_0}(t)$ to denote the integral curve as in (\ref{est-int-curve}), but with $\wh f$ replaced by the bias-corrected version $\wh f^{\,bc}$.

\subsubsection{Confidence regions for ${\cal M}$ and ${\cal L}$ using local adjustment by the gradient}

Using the de-biased estimator $\wh f^{\,bc}$ along with the integral curve $\wh\X^{bc}_{x}(t)$, we now present the  construction of two confidence regions, one based on asymptotic distribution theory, and the other based on the bootstrap. The latter has a faster rate of approximation of the coverage probability. \\[5pt]
For $x \in \R^d,$ define $\wh \theta^{\,bc}_{x} \in \R$ through the property
\begin{align}
\wh f^{\,bc} (\wh \X^{bc}_{x} ( \wh\theta^{\,bc}_{x} ) )=c.\label{thetatilde}
\end{align}
For large sample size the existence and uniqueness of $\wh \theta^{\,bc}_{x}$ for $x\in\mathcal{M}$ are proved in Lemma~\ref{existenceuniqueness}. For finite sample, in case the solution to (\ref{thetatilde}) is not unique, we take $\wh\theta^{\,bc}_{x}$ as the infimum of the set of solutions; and whenever there is no solution to (\ref{thetatilde}), we set $\wh\theta^{\,bc}_{x}$ to be the smallest value of $\argmin_{\theta} |\wh f^{\,bc} (\wh \X^{bc}_{x} ( \theta ) )-c|$. Noticing that $\wh{\X}^{bc}_{x} ( \wh\theta^{\,bc}_{x} ) \in \wh {\cal M}^{\,bc} = (\wh f^{\,bc})^{-1}(c)$, we now set $\wh x$ in (\ref{vert-horiz}) as $\wh x = \wh\X^{bc}_{x}(\wh\theta^{\,bc}_{x})$. \\[5pt]
Letting $0 < \alpha < 1,$ and recalling the definition of $\wh a^{(d)}_{1- \alpha}$ given in (\ref{def-hata}), we define
\begin{align*}
\wh C_{n,4}(1-\alpha) = \big\{x \in \R^d:\, \|\nabla \wh f^{\,bc}(\wh \X^{bc}_{x}(\wh \theta^{\,bc}_x))\|\,\|\wh \X^{bc}_{x}(\wh \theta^{\,bc}_x) - x\| \le  \wh a^{(d)}_{1- \alpha} \big\}.
\end{align*}
It can be shown (see Lemma~\ref{existenceuniqueness}) that $x\mapsto \wh \X^{bc}_{x}(\wh \theta^{\,bc}_x)$ defines a bijective map between $\mathcal{M}$ and $\wh{\mathcal{M}}^{bc}$ when the sample size is large enough. Thus, for large sample size, we can equivalently write  
$$\wh C_{n,4}(1-\alpha) = \{\wh \X^{bc}_{z}(t):\, \|\nabla \wh f^{\,bc}(z)\|\,\|z - \wh \X^{bc}_{z}(t) \| \le  \wh a^{(d)}_{1- \alpha}, z \in \wh {\cal M}^{\,bc},\,t \in \R \}.$$  
This also indicates an algorithm: For a dense enough subset of values $z \in \wh {\cal M}^{\,bc}$, run the integral curve $\wh \X^{\,bc}_{z}(t)$, and check whether the condition in the definition of the confidence region holds. A bootstrap version of this confidence region is given by
\begin{align*}
\wh C^*_{n,4}(1-\alpha) = \big\{x \in \R^d:\, \|\nabla \wh f^{\,bc}(\wh \X^{bc}_{x}(\wh \theta^{\,bc}_x))\|\,\|\wh \X^{bc}_{x}(\wh \theta^{\,bc}_x) - x\| \le  \wh c^{*,E}_{1- \alpha} \big\},
\end{align*}
where $ \wh c^{*,E}_{1- \alpha} $ is as in (\ref{Cn2}). \\

Let $\wh {\mathcal{L}}^{bc} =\{x\in\mathbb{R}^d: \; \wh f^{bc}(x) \geq c\}$ and define
\begin{align*}
&\wh C_{n,4}^{-}(1-\alpha) = \wh {\mathcal{L}}^{bc} \cup \wh C_{n,4}(1-\alpha), \quad\quad \wh C_{n,4}^{+}(1-\alpha) = \wh {\mathcal{L}}^{bc} \backslash \wh C_{n,4}(1-\alpha) ,\\
&\wh C_{n,4}^{*,-}(1-\alpha) = \wh {\mathcal{L}}^{bc} \cup \wh C_{n,4}^*(1-\alpha) \quad\quad \wh C_{n,4}^{*,+}(1-\alpha) = \wh {\mathcal{L}}^{bc} \backslash \wh C_{n,4}^*(1-\alpha) .
\end{align*}

\begin{theorem}\label{asympt-Cn3-Cn4}
{\bf Part 1.} Let $0 < \alpha < 1.$ Suppose that ({\textbf F1}), ({\textbf F2}), ({\textbf K}), ({\textbf A}), ({\textbf H1})$^2$ and ({\textbf H2})$^2$ hold. Then we have, as $n \to \infty$,
\begin{align}
&\mathbb{P}\Big( \mathcal{M} \subset \wh C_{n,4}(1-\alpha)  \Big) =1-\alpha + o(1),\quad \text{and}\label{cn4confm}\\
&\mathbb{P}\Big( \wh C_{n,4}^+(1-\alpha) \subset \mathcal{L} \subset \wh C_{n,4}^-(1-\alpha)  \Big) =1-\alpha + o(1).\label{cn4confl}
\end{align}
{\bf Part 2.} Let $0 < \alpha < 1,$ and suppose that (\textbf{F1}), ({\textbf F2}), (\textbf{K}), ({\textbf H1})$^2$ and ({\textbf H2})$^2$ hold. Let $\gamma_n^{bc}$ as in (\ref{gammanbc}). If $\Psi_n(\gamma_n^{bc}) = o(1)$ as $n\rightarrow\infty$, then we have
\begin{align}
& \mathbb{P} \left( \mathcal{M} \subset \wh C_{n,4}^*(1-\alpha) \right) = 1-\alpha + O\big(\Psi(\gamma_n^{bc})\big),\quad \text{and} \label{cn4starconfm}\\
& \mathbb{P} \left( \wh C_{n,4}^{*,+}(1-\alpha) \subset \mathcal{L} \subset \wh C_{n,4}^{*,-}(1-\alpha) \right) = 1-\alpha + O\big(\Psi(\gamma_n^{bc})\big).\label{cn4starconfl}
\end{align}
\end{theorem}
Notice that we do not have rate of convergence for the coverage probability of the confidence regions based on the large sample theory. At this point it is an open question whether, and if yes, how, to derive such rates of approximations. Since we use the approximation of the extreme value distribution of Gaussian fields, it is expected to be a slow rate $(\log{n})^{-1}$, following a similar argument given in Hall (1979).

\subsubsection{Confidence regions for ${\cal M}^E$ and ${\cal L}^E$ adjusted by gradient}\label{meleabg}
Similar to the above, confidence regions for ${\cal{M}}^E$ and ${\cal{L}}^E$ can be constructed and analysed. For instance, a bootstrap confidence region for ${\cal{M}}^E$ is given by 
\begin{align*}
\wh C^{*,E}_{n,4}(1-\alpha) = \big\{x \in \R^d:\, \|\nabla \wh f(\wh \X_{x}(\wh \theta_x))\|\,\|\wh \X_{x}(\wh \theta_x) - x\| \le  \wh c^{*,E}_{1- \alpha} \big\}.
\end{align*}
The lower and upper ``bounds'' of the confidence region for ${\cal L}^E$ can be constructed by taking set difference and union between $\wh{\cal L}$ and $\wh{C}^{*,E}_{n,4}(1-\alpha)$, respectively. Note that the construction uses estimators based on $\wh f$, instead of $\wh f^{\,bc}$, similar to $\wh{C}^{*,E}_{n,2}(1-\alpha)$ in (\ref{Cn2}). It can be shown that their rates of convergence for the coverage probability are $O\big(\Psi(\gamma_n^{E})\big)$ if we assume that (\textbf{F1}), ({\textbf F2}), (\textbf{K}), ({\textbf H1})$^2$ and ({\textbf H2})$^2$ hold and $\Psi(\gamma_n^{E})) =o(1)$. The proof follows the same arguments given in the proof of Theorem~\ref{asympt-Cn3-Cn4}. Details are omitted. 

\subsubsection{Confidence regions not locally adjusted by the gradient}
The confidence regions constructed here are related to the confidence regions constructed in the previous subsection, but in contrast to them, here the dependence on the estimated gradient is more indirect through the construction of the integral curve. As a result, the width of the confidence regions only depends on the length of the integral curve, and it is not locally adjusted by the gradient. The construction is as follows. Recall that $\wh f^{*,E}(x) = {\mathbb E}^* \wh f^*(x),$ and let 
\begin{align*}  
\wh{\mathcal{M}}^{*,E} = \{x:\;  \wh f^{*,E}(x)=c\}.
\end{align*}
Let $\wh\X_x^{*}$ be the integral curve driven by $\nabla \wh f^{*}$, and let $\widehat{\theta}_x^{*}$ be the first time $t$ at which $\widehat{\X}_x^{*}(t)$ hits $\mathcal{M}^{*}$. In the case of large sample size the existence and uniqueness of $\wh\theta_x^{*}$ follow the same argument as for $\wh \theta_x$. 
For $0 < \alpha < 1,$ let $\wh d_{1-\alpha}^{*,E}$ be the quantile of order $(1-\alpha)$ 
for $\sup_{x\in\wh{\mathcal M}^{*,E}}|\wh \X_x^{*}(\wh\theta_x^{*})-x|,$ and define the set
$$\wh C^*_{n,5}(1-\alpha) = \{\wh \X^{bc}_x(t):\;\|\wh \X^{bc}_x(t)-x\|\leq\wh d_{1-\alpha}^{*,E},\; x\in\wh{\mathcal M}^{\,bc}, t \in \R\}.$$
Also define
\begin{align*}
\wh C_{n,5}^{*,-}(1-\alpha) = \wh {\mathcal{L}}^{bc} \cup \wh C_{n,5}^*(1-\alpha) \quad\text{and}\quad \wh C_{n,5}^{*,+}(1-\alpha) = \wh {\mathcal{L}}^{bc} \backslash \wh C_{n,5}^*(1-\alpha) .
\end{align*}

Let $\zeta_n = \beta_{n,g}^{(1)}\log{n} + \beta_{n,g}^{(0)}\left[h^{-1}\log{n}+ \sqrt{nh^{d+4}\log{n}}\right]$. Now we have the following result:

\begin{theorem}\label{asympt-Cn5}
Suppose that (\textbf{F1}), (\textbf{F2}), (\textbf{K}), ({\textbf H1})$^2$ and ({\textbf H2})$^2$ hold. Then, if $\Psi_n(\gamma_n^{bc}) +\zeta_n= o(1)$ as $n\rightarrow\infty,$ where $\gamma_n^{bc}$ is as in (\ref{gammanbc}), we have as $n\rightarrow\infty$,
\begin{align}
&\mathbb{P} \left( \mathcal{M} \subset \wh C_{n,5}^*(1-\alpha) \right) = (1-\alpha) + O\left(\Psi_n(\gamma_n^{bc}) + \zeta_n\right), \quad \text{and}\label{cn5confm}\\
& \mathbb{P} \left( \wh C_{n,5}^{*,+}(1-\alpha) \subset \mathcal{L} \subset \wh C_{n,5}^{*,-}(1-\alpha) \right) = 1-\alpha + O\left(\Psi_n(\gamma_n^{bc}) + \zeta_n\right).\label{cn5confl}
\end{align}
\end{theorem}
\begin{remark}{\em
Note that if we further assume that $g$ and $h$ are of the same rate, i.e., there exist $0<C_1,C_2<\infty$ such that $C_1<h/g<C_2$ as $n\rightarrow\infty$, then $\zeta_n$ can be absorbed into $\Psi_n(\gamma_n^{bc})$ in the above results.}
\end{remark}
Similar to subsection~\ref{meleabg}, a confidence region for ${\cal M}^E$ can be constructed by
$$\wh C^{*,E}_{n,5}(1-\alpha) = \{\wh \X_x(t):\;\|\wh \X_x(t)-x\|\leq\wh d_{1-\alpha}^{*,E},\; x\in\wh{\mathcal M}, t \in \R\},$$
and the confidence region for ${\cal L}^E$ has `upper and lower boundaries' given by
\begin{align*}
\wh C_{n,5}^{*,E,-}(1-\alpha) = \wh {\mathcal{L}} \cup \wh C_{n,5}^{*,E}(1-\alpha) \quad\text{and}\quad \wh C_{n,5}^{*,E,+}(1-\alpha) = \wh {\mathcal{L}} \backslash \wh C_{n,5}^{*,E}(1-\alpha) .
\end{align*}

\subsection{Horizontal variation based methods not based on integral curves}\label{Horizontal2}
It follows from the Tubular Neighborhood Theorem (e.g. see Theorem 11.4 of Bredon, 1993) that for all $x\in\mathcal{M}$ and $n$ large enough, there exist unique $X_x\in\wh{\mathcal{M}}^{bc}$ and $s_{x}\in\mathbb{R}$ such that $x=X_x+s_{x}\nabla \widehat f^{bc}(X_x)$. Similarly, for all $x\in\wh{\mathcal{M}}^{*,E}$ and large sample size, we can find $X_x^*\in\wh{\mathcal{M}}^*$ and $s_x^*$ be such that $x=X_x^* + s_x^*\nabla \widehat f^*(X_x^*)$. \\

Now let $\widehat b_{1-\alpha}^{*,E}$ be the quantile of $(1-\alpha)$ for $\sup_{x\in\wh{\mathcal{M}}^{*,E}}\|X_{x}^*-x\|$, and for $0 < \alpha < 1$ define
$$ \wh C^*_{n,6}(1-\alpha) = \left\{x+t\nabla \widehat f^{bc}(x): \; \|t\nabla \widehat f^{bc}(x)\|\leq \widehat b_{1-\alpha}^{*,E}, \;\; x\in\wh{\mathcal{M}}^{bc}\right\}.$$
Also define
\begin{align*}
\wh C_{n,6}^{*,-}(1-\alpha) = \wh {\mathcal{L}}^{bc} \cup \wh C_{n,6}^*(1-\alpha) \quad\text{and}\quad \wh C_{n,6}^{*,+}(1-\alpha) = \wh {\mathcal{L}}^{bc} \backslash \wh C_{n,6}^*(1-\alpha) .
\end{align*}
Under suitable regularity conditions $\wh C^*_{n,6}(1-\alpha)$ is a confidence region for ${\cal M}$ of asymptotic coverage level $1-\alpha$, and $\wh C_{n,6}^{*,-}(1-\alpha)$ and $\wh C_{n,6}^{*,+}(1-\alpha)$ give the upper and lower ``bounds'' of an asymptotic $(1-\alpha)$ confidence region for $\mathcal{L}$. The convergence rate of the coverage probability of $\wh C^*_{n,6}(1-\alpha)$ for $\mathcal{M}$ as well as $\{\wh C_{n,6}^{*,+}(1-\alpha), \wh C_{n,6}^{*,-}(1-\alpha)\}$ for $\mathcal{L}$ can be derived in the way similar to the proof of Theorem~\ref{asympt-Cn5}. However, we do not further pursue it here. 
While the geometric construction of the confidence region $\wh C^*_{n,6}(1-\alpha)$ is essentially the same as the one constructed in Chen et al. (2017), our derivation via the tubular neighborhood theorem provides a slightly different angle to the construction. 
A similar confidence region for ${\cal M}^E$ is given by 
$$ \wh C^{*,E}_{n,6}(1-\alpha) = \left\{x+t\nabla \widehat f(x): \; \|t\nabla \widehat f(x)\|\leq \widehat b_{1-\alpha}^{*,E}, \;\; x\in\wh{\mathcal{M}}\right\}.$$
For ${\cal L}^E$, we can use the following upper and lower ``bounds'':
\begin{align*}
\wh C_{n,6}^{*,E,-}(1-\alpha) = \wh {\mathcal{L}} \cup \wh C_{n,6}^{*,E}(1-\alpha) \quad\text{and}\quad \wh C_{n,6}^{*,E,+}(1-\alpha) = \wh {\mathcal{L}} \backslash \wh C_{n,6}^{*,E}(1-\alpha).
\end{align*}
Rates of convergence of the coverage probabilities for these confidence regions can be derived by using similar ideas as above. No further details are given.

\section{Performance of confidence regions and geometry}\label{Discussion}

The above discusses the large sample behavior of various confidence regions for level sets. Bootstrap based methods show a faster rate of convergence of their coverage probability to the nominal level than the methods based on asymptotic distribution theory, which is not a surprise. A more detailed comparison based on the theoretical developments is not entirely straightforward, because different confidence sets depend on a different number of bandwidths to be chosen, requiring different assumptions, etc.  However, for finite samples, certain relations between the geometry of the underlying density and the performance of the different types of confidence regions give some interesting insight. This will be discussed now. \\[5pt]
1. By construction, most of the confidence regions for ${\cal M}$ or ${\cal M}^E$ based on horizontal variation constitute a band (or a tube) of constant width about the estimated target isosurface. The width of the tube depends on the global behavior of the density in a neighborhood about the targeted isosurface. Since horizontal variation based methods are essentially based on the worst case behavior (similar to the supremum distance), it can be expected that for densities for which the norm of the gradient varies a lot along the isofurface, this confidence band tends to be unnecessarily wide. This is also illustrated in our simulation study in Section~\ref{simulations}. 

In contrast to that behavior of the horizontal confidence regions, the width of the vertical variation based confidence regions has local adaptivity. Essentially, their width at a given point of the isosurface is inversely proportional to the norm of the gradient, and thus they are containing additional information about the geometry of the density. \\[5pt]
2. Our theoretical assumptions restrict the level $c$ to be strictly larger than zero, which means that in local neighborhoods about the corresponding isosurface the density is bounded away from zero. However, a relatively small level $c$ can still provide problems in finite samples. For instance, vertical variation based methods of the form $\wh f[c - \wh a_n, c + \wh a_n]$ (or the ones using the bias-corrected density estimator) might be very large in volume, as the lower bound $c - \wh a_n$ might be less than zero, meaning that the outer confidence region only ends at the support of the density estimator used. Nevertheless, the probability content carried by the confidence regions might still be small. (Asymptotically, this problem of course disappears simply because $\wh a_n$ converges to zero.) In such situations, the volume of horizontal variation based methods tend to be of smaller volume than the ones of the vertical based methods, but the probability mass carried by them might nevertheless be larger. This can be seen in the simulation results presented in Table~\ref{simulationres} when inspecting the column corresponding to Case 3. \\[5pt]
3. Another interesting scenario corresponds to levels $c$ that are close to critical values of the density. Similar to the previous item, this problem does not appear in a large sample scenario, because our assumptions require the gradient along the iso-surface to be bounded away from zero. For finite samples, however, we observe the following interesting geometric challenge. 
\begin{figure}[h]
\centering
\includegraphics[height=3in,width=3.5in]{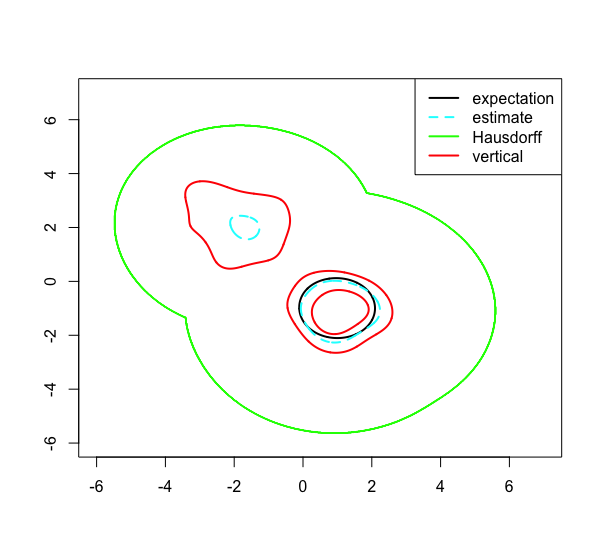}
\vspace*{-0.3cm}
\caption{An example of Case 4 (given below); $n=200$. The target is $\mathcal{M}^E$ (black curve). Notice, $\mathcal{M}^E$ only has one connected component, while $\mathcal{M}$ has two (not shown here). The estimate $\wh{\mathcal{M}}$ (cyan dotted curve) has two connected components. Red curves give the boundaries of the 90\% confidence region using the vertical method corresponding to V.e below. Green curve are the boundary of the 90\% confidence region using the Hausdorff method.}
\label{Fig:compar}
\end{figure}

Suppose that $c$ is only slightly smaller than a level at which the true density has a local maximum that is not a global maximum, and let $x_0$ denote the point at which this local maximum is attained. Then, under our regularity assumptions, there will be a neighborhood of $x_0$ that is part of the superlevel set ${\cal L}$. However, the value of the density estimates $\wh f(x_0)$ or $\wh f^{\,bc}(x_0)$ might, due to random fluctuation, not exceed the value $c$, so that the estimated superlevel set does not contain a neighborhood about $x_0$. Since the confidence regions based on horizontal variation are built around the contours of these estimated superlevel sets, they might miss such areas - and note that these areas then could in fact lie far from the confidence region. In such cases the missed target regions tend to be small in size. Nevertheless, the topology (homology) of the confidence regions will in general be different from the one of the target contour, and the normal compatibility assumption (e.g. see Chazal et al. 2007) used in Chen et al. (2017) for the construction of horizontal variation based confidence regions will be violated in such cases.

Horizontal variation based confidence regions will not perform well in this scenario. Observe that the quantiles used in their construction are essentially based on the maximum distance of the estimated and the true contours. This distance tending to be large means that these quantiles will tend to become large, leading to wide confidence tubes. By contrast, the vertical distance is less impacted by the different topology of the estimated and true contours, and hence the vertical variation based confidence regions suffer less from having very large volume in this scenario.

\begin{figure}[ht]
\centering
\includegraphics[width=4.8in, height=2.6in]{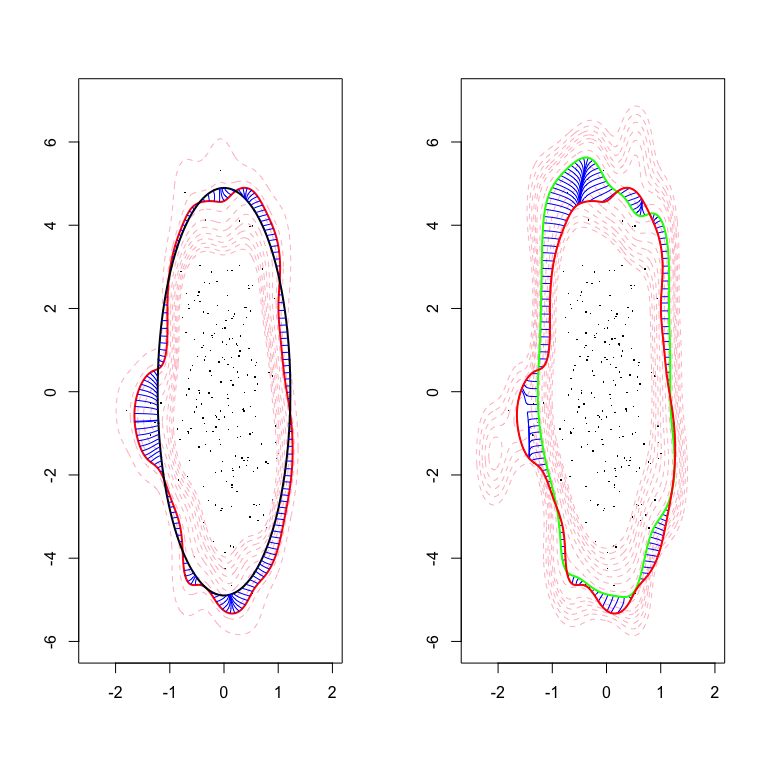}
\caption{Challenges of constructing confidence regions based on integral curves are illustrated.}
\label{Fig:horiz}
\end{figure}

A scenario as the one just discussed is included in the simulation study presented below. See Table 1, Case 4, where one can see that the Hausdorff based methods tend to be quite large. A similar remark applies to the integral curve based methods, such as $\wh C_{n,4}(1-\alpha)$, where the indicated problem is expressed by the non-existence of $\wh\theta_x^{\,bc}$ for a non-negligible set of starting values $x$. A similar discussion applies when the level $c$ is narrowly above a critical level, illustrated in Figure~\ref{Fig:compar}.\\[5pt]
4. Confidence regions based on estimated integral curves might, for finite sample size, suffer from the local geometry of the kernel density estimator, as is illustrated in Figure~\ref{Fig:horiz}. The two panels in this figure show a scenario (for finite sample size) that violates the bijective condition for horizontal methods, which is shown to hold for large samples (see Lemma~\ref{existenceuniqueness}). 
The violation is due to small local minimum of the kernel density estimator, and this local geometric property then `diverts' the integral curves from their expected path. This also gives rise to numerical challenges. In the situation shown in the figure, a sample of size 200 was drawn from a density function in (\ref{simdens}) below with $a=2,$ and then a bootstrap sample was drawn. The focus was on the level set corresponding to $p=0.95$. In the left panel, the pink curves are contour lines of $\wh f$; the red curve is $\wh{\mathcal{M}}$; the black curve is $\mathcal{M}$; the blue curves are trajectories of integral curves driven by the gradients of $\wh f$. In the right panel, the pink curves are contour lines of $\wh f^*$; the red curve is still $\wh{\mathcal{M}}$; the green curve is $\mathcal{M}^*$; the blue curves are trajectories of integral curves driven by the gradients of $\wh f^*$. Notice that the trajectories fail to define a bijective mapping between $\wh{\mathcal{M}}$ and $\wh{\mathcal{M}}^*$ around (-1.5, -1.5) due to the existence of a local minimum of $\wh f^*$.

\section{Simulations}\label{simulations}
This simulation study compares six (one large sample based, and five bootstrap based) confidence regions for ${\cal M}$ and ${\cal M}^E$ in terms of coverage probability and volume, thereby
%
\begin{itemize}
\item considering different levels $c$ (low, high, close to critical levels), and
\item comparing vertical variation based and horizontal variation based methods.
\end{itemize}
%
All the horizontal variation based methods are more computationally involved than the vertical methods. A preliminary simulation study showed our horizontal variation based methods to behave similarly to the Hausdorff-distance based approach of Chen et al. (2017). Therefore we here only use the latter to represent the horizontal variation based methods. With $d(x, \wh {\cal M}) = \inf_{s \in \wh {\cal M}} \|x - s\|,$ these confidence regions have the form $\big\{ x\in \R^d: d(x, \wh {\cal M}) \le \wh e^*_{1-\alpha}\big\},$ where $\wh e^*_{1-\alpha}$ is a bootstrap based estimate of the $(1-\alpha)$-quantile of $d(x, \wh {\cal M}).$ Recall that the Hausdorff distance between ${\cal M}$ and $\wh {\cal M}$ is given by 
$$d_H({\cal M}, \wh{\cal M}) = \max\Big( \sup_{x \in {\cal M}}d(x, \wh {\cal M}),\,\sup_{s \in \wh {\cal M}}d(s, {\cal M})\Big), $$
and thus, $\sup_{x \in {\cal M}}d(x, \wh {\cal M})$ is `one part' of the Hausdorff distance. However, if ${\cal M}$ and ${\wh {\cal M}}$ are {\em normal compatible}, then $\sup_{x \in {\cal M}}d(x, \wh {\cal M}) = \sup_{s \in \wh {\cal M}}d(s, {\cal M})$, and 
$$d_H({\cal M}, \wh{\cal M}) = \sup_{x \in {\cal M}}d(x, \wh {\cal M}).$$
Chen et al. (2017) are using this approach  with ${\cal M}$ replaced by ${\cal M}^E$, and they show that, under certain regularity assumptions, normal compatibility of ${\cal M}^E$ and ${\wh {\cal M}}$ holds asymptotically with probability tending to one. \\

One of the models used in our simulations is the bivariate normal with 
\begin{align}\label{simdens}
f(x,y;a) = \frac{1}{2\pi}e^{\frac{-a^2x^2-y^2/a^2}{2}},
\end{align}
where the contours of the density function are ellipses with $a$ defining their eccentricity.  
Then, for $c=f(x_0,y_0;a)$ with $a^2x_0^2+y_0^2/a^2=r_0^2$ for some $0<r_0<\infty$, the probability over the superlevel set $\{(x,y):\;\;f(x,y;a)\geq c\}$ is
\begin{align*}
p=\int_{a^2x^2+y^2/a^2\leq r_0^2} f(x,y)dxdy =\int_{0}^{2\pi}\int_0^{r_0}\frac{1}{2\pi}e^{-r^2/2}rdrd\theta =1-2\pi c.
\end{align*}
We choose $p=50\%$ and $p=95\%$. Our second model is a mixture of normal distributions of the form
\begin{align}\label{simmix}
0.5\mathcal{N}\left((-2,2)^T,\;1.5\mathbf{I}_2\right) \; + \; 0.5\mathcal{N}\left((1,-1)^T,\;0.5\mathbf{I}_2\right)
\end{align}
This density has two modes with corresponding heights 0.065 and 0.11, respectively. In our study we chose the level $c=0.048$, which lies slightly below the lower local maximum of the mixture of normals (cf. 3. in Section~\ref{Discussion}). We consider the following 4 cases:\\
\indent Case 1: density in (\ref{simdens}) with $a=1$ and $p=0.5$,\\
\indent Case 2: density in (\ref{simdens}) with $a=2$ and $p=0.5$,\\
\indent Case 3: density in (\ref{simdens}) with $a=1$ and $p=0.95$,\\
\indent Case 4: density in (\ref{simmix}) with $c=0.048$.\\

 As a kernel we choose the form
\begin{align*}
K(x,y) = \left(\frac{693}{512}\right)^2(1-x^2)^5(1-y^2)^5\mathbf{1}{\{|x| \leq 1 \text{ and } |y|\leq 1\}}.
\end{align*}
We ran the simulation for 400 times. In each iteration, a sample of size $n$ was randomly drawn from the given distribution and then a bootstrap procedure based on 250 bootstrap re-samplings was performed to create the confidence regions using the following methods:\\[5pt]
 (H)\;Hausdorff-distance based approach of Chen et al. (2017) for the smoothed level set;\\[3pt]
(V.e)\; vertical variation based confidence region $\wh C_{n,2}^{*,E}(1-\alpha)$ for the smoothed level set;\\[3pt]
(V)\; vertical variation based confidence region $\wh C_{n,2}^*(1-\alpha)$ for the true level set;\\[3pt]
(V.bc)\; vertical variation based confidence region with bias correction $\wh C_{n,3}^*(1-\alpha)$ for the true level set;\\[3pt]
(V.us)\; vertical variation based confidence region with undersmoothing $\wh C_{n,2}^{*,E}(1-\alpha)$ for the true level set (see Remark~\ref{vert-boot});\\[3pt]
(V.ls)\; vertical variation based large sample confidence region $\wh C_{n,1}(1-\alpha)$ for the true level set.\\[5pt]
The bandwidths involved in the construction of these confidence regions are selected using the direct plug-in method. In particular, we use the plug-in optimal bandwidth for kernel density estimation as $h$ and $g$, while using the plug-in optimal bandwidth for the second derivative estimation as $l$. In fact, we choose different bandwidths for each of the two dimensions. For (V.us), we used 70\% of the optimal plug-in bandwidths.\\

The confidence level was set to be 90\%. With the 400 runs, we calculated the coverage probabilities (C.P.) of these confidence regions as well as their average Lebesgue measures ($\bar\lambda$) and average probability measures ($\bar P$). Note that in the general form $\wh f[c-\wh a, c+ \wh a]$ or $\wh f^{bc}[c-\wh a, c+ \wh a]$ of the confidence regions based on vertical variation, sometimes $c-\wh a<0$ for Case 3. Since the kernel function $K$ we used has bounded support, so do $\wh f$ and $\wh f^{bc}$. The outer boundary of the confidence regions based on vertical variation for the level sets is in fact the support of the density estimator (cf. 2. of Section~\ref{Discussion}). For numerical reasons, we used $\wh f[\max(c-\wh a,\omega), c+ \wh a]$ or $\wh f^{bc}[\max(c-\wh a,\omega), c+ \wh a]$ as the confidence regions, where we took $\omega=10^{-6}$.\\

\begin{table}[h]
\centering
\caption{Simulation results}
\label{simulationres}
\vspace*{0.2cm}
\resizebox{4.8in}{!}{
\begin{tabular}{ll|lll|lll|lll|lll}
\hline\hline
                        &             & \multicolumn{3}{c|}{Case 1} & \multicolumn{3}{c|}{Case 2} & \multicolumn{3}{c|}{Case 3} & \multicolumn{3}{c}{Case 4} \\
                        &             & C.P.      & $\bar\lambda$    & $\bar P$      & C.P.      & $\bar\lambda$     & $\bar P$      & C.P.      & $\bar\lambda$     & $\bar P$      & C.P.      & $\bar\lambda$    & $\bar P$      \\ \hline\hline
\multirow{4}{*}{n=200}  
& H	&0.98	&8.40	&0.63	&1.00	&13.18	&0.82	&1.00	&41.21	&0.49	&1.00	&59.46	&0.88\\
& V.e	&0.87	&6.91	&0.53	&0.91	&6.87	&0.53	&0.95	&35.16	&0.15	&0.90	&12.17	&0.47\\
& V	&0.91	&7.33	&0.56	&0.94	&7.28	&0.56	&0.95	&36.17	&0.17	&0.85	&13.14	&0.49\\
& V.bc	&0.50	&6.28	&0.46	&0.56	&6.26	&0.46	&0.66	&20.63	&0.15	&0.53	&11.75	&0.43\\
& V.us	&0.94	&13.55	&0.77	&0.96	&13.45	&0.77	&0.96	&28.58	&0.26	&0.88	&29.82	&0.71\\\hline\hline
\multirow{4}{*}{n=1000} 
& H	&0.95	&4.26	&0.34	&0.98	&7.75	&0.6	&1.00	&21.22	&0.18	&0.98	&33.35	&0.76\\
& V.e	&0.90	&3.77	&0.30	&0.89	&3.79	&0.30	&0.92	&36.31	&0.09	&0.9	&7.41	&0.31\\
& V	&0.94	&4.09	&0.33	&0.94	&4.11	&0.33	&0.92	&37.37	&0.11	&0.89	&8.00	&0.33\\
& V.bc	&0.61	&3.54	&0.27	&0.61	&3.55	&0.28	&0.64	&25.60	&0.1	&0.63	&7.41	&0.31\\
& V.us	&0.94	&6.27	&0.47	&0.93	&6.30	&0.47	&0.96	&30.81	&0.15	&0.94	&13.78	&0.49\\\hline\hline
\multirow{4}{*}{n=5000} 
& H	&0.93	&2.32	&0.19	&0.97	&4.47	&0.36	&0.98	&9.67	&0.07	&0.99	&15.29	&0.58\\
& V.e	&0.90	&2.17	&0.17	&0.88	&2.16	&0.17	&0.94	&9.38	&0.06	&0.90	&4.77	&0.22\\
& V	&0.94	&2.38	&0.19	&0.96	&2.37	&0.19	&0.96	&13.51	&0.07	&0.92	&5.17	&0.23\\
& V.bc	&0.74	&2.08	&0.16	&0.73	&2.07	&0.16	&0.74	&8.82	&0.06	&0.71	&4.84	&0.22\\
& V.us	&0.97	&3.52	&0.28	&0.96	&3.49	&0.27	&0.97	&33.49	&0.10	&0.95	&7.89	&0.33\\
& V.ls &1.00	&5.96	&0.44	&1.00	&6.39	&0.47	&1.00	&34.45		&0.13	&1.00	&13.85	&0.50\\\hline\hline
 \multirow{4}{*}{n=50000} & H  &0.93		&1.04	&0.08	&0.92	&2.06	&0.17	&0.94	&3.86	&0.03	&1.00	&7.89	&0.36\\
			& V.e	&0.92	&1.01	&0.08	&0.92	&0.98	&0.08	&0.90	&3.68	&0.03	&0.93	&2.73	&0.13\\
			& V		&0.96	&1.11	&0.09	&0.95	&1.09	&0.09	&0.94	&4.65	&0.03	&0.97	&2.95	&0.14\\
			& V.bc	&0.81	&0.98	&0.08	&0.84	&0.96	&0.08	&0.77	&3.57	&0.03	&0.83	&2.75	&0.13\\
			& V.us	&0.97	&1.63	&0.13	&0.96	&1.60	&0.13	&0.96	&6.62	&0.05	&0.98	&4.09	&0.19  \\ 
			&V.ls		&1.00	&1.64	&0.13	&1.00	&1.72	&0.14	&1.00	&6.44	&0.05	&1.00	&4.09	&0.19\\\hline    \hline             
\end{tabular}
}
\end{table}

In Table~\ref{simulationres}, (H) has to be compared with (V.e) since these methods are both targeting the smoothed level sets $\mathcal{M}^E$. Overall it is clear that the vertical method (V.e) outperforms the horizontal method (H). Detailed discussions have been given in Section~\ref{Discussion}.\\

We only include results for the large sample confidence regions (V.ls) with $n=5000$ and $n=50000$. This is because the formula in Theorem 2.1 requires $h<1$, which cannot be satisfied when $n$ is small in our examples. Overall the large sample confidence regions have conservative coverage probabilities in our examples. This is not surprising because it is well-known that the convergence rate of the coverage probability is slow for such large sample confidence regions. However, their volumes are comparable to those of the bootstrap confidence regions when the sample size is large, which makes (V.ls) a competitive option considering its computation does not require bootstrap.\\

Comparing (V), (V.bc) and (V.us), which are all confidence regions for the true level sets $\mathcal{M}$, it is apparent that (V) performs best. Convergence of the coverage probability of bias correction method (V.bc) is the slowest. This seems to be caused by slow convergence in the estimation of the second order derivatives in the bias correction. The confidence regions based on the undersmoothing method (V.us) have the largest volume. This is because the variance becomes large when the selected bandwidth is small.

\section{Conclusion} \label{Conclusion} We have constructed and analyzed various confidence regions for density superlevel sets and density isosurfaces based on plug-in estimates using kernel density estimation. The analysis is done in terms of large sample theory and also in the finite sample setting using simulations. The geometry underlying the construction of the different types of confidence regions is discussed. Geometric considerations also play a role in the interpretation of the finite sample behavior of the confidence regions. Overall, vertical variation based confidence intervals appear to have an edge over the horizontal methods. \\[5pt]
The kernel estimator used in our investigations can of course be replaced by other (non-parametric) density estimators. For such modifications, the large sample behavior of the corresponding coverage probabilities might be analyzed using a similar Ansatz as in this work (see discussion of ``Structure of the rates of convergence for the coverage probabilities of bootstrap based confidence sets'' given in Section 2.2). This requires the investigation of all the relevant properties needed for our approach to go through.  \\[5pt]
There are various open questions related to the construction of confidence regions for density level sets. For instance, what can be said about optimality of the rates of convergence of the coverage probabilities? (Thanks to the referee for asking this question). Besides some classical work by Hall and Jing (1995), the only other work related to this question we are aware of is Calonico et al. (2018b). The role of bias correction in this context might be explored as well. For some recent work on bias correction see Chen (2017), and Calonico et al. (2018a). While we have been concentrating on density level sets, a similar approach might work for level sets of other functions, such as regression level sets, for instance. Level sets also play an integral role in the context of topological data analysis (persistent homology). Similar to our vertical variation based upper and lower confidence sets, Bobrowski et al. (2017) are using such upper and lower approximations of the level set to construct estimates for the topology of a single density level set. It might be worthwhile to explore this connection in more detail.

\section{Proofs}\label{Proofs}

\subsection{Proof of Theorem~\ref{LevelSetConf}}

A key ingredient to the proof of Theorem~\ref{LevelSetConf}, is the following special case of the main theorem in Qiao and Polonik (2018):


\begin{theorem}\label{ProbabilityResult}
Let $\mathcal{H}\subset\mathbb{R}^d$ ($d \geq 2$) be a compact set. Let $Z_h(x), x\in \mathcal{H}$, $0<h\leq1$ be a sequence of centered Gaussian fields with covariance
\begin{align}\label{varfunc}
r_h(x+\Delta x, x) = 1 - h^{-2}\|D \Delta x\|^2 + o(h^{-2}\|\Delta x\|^2),
\end{align}
uniformly in $h\in(0,1]$ and $x\in\mathcal{H}$ as $\Delta x/h \rightarrow 0$, where $D$ is $d\times d$ positive definite matrix. Let $r<d$ and $\mathcal {M}\subset\mathcal{H}$ be a $r$-dimensional compact Riemannian manifold with reach $\Delta(\mathcal{M})>0$. For any $\delta>0$, define
\begin{align*}
Q(\delta):=\sup_{0<h\leq 1}\{|r_h(x+\Delta x,x)|: x+\Delta x\in\mathcal {M}, x\in\mathcal {M}, \|\Delta x\|>h\delta\}.
\end{align*}
Suppose for any $\delta>0$, there exists a positive number $\eta$ such that
\begin{align}\label{SupGauss1}
Q(\delta) < \eta<1,
\end{align}
In addition, assume that there exists $\eta>0$ and $\delta_0$, such that, for any $\delta>\delta_0$, we have 
\begin{align}\label{SupGauss2}
Q(\delta) |(\log{\delta})^r| \leq (\log{\delta})^{-\eta}.
\end{align}

For any fixed $z$, define
\begin{align}\label{ThetaExp}
\phi(z) &=\sqrt{2r\log{h^{-1}}}+\frac{1}{\sqrt{2r\log{h^{-1}}}}\bigg[z+\Big(\frac{r}{2}-\frac{1}{2}\Big)\log{\log{h^{-1}}}\nonumber\\
&\hspace{3cm}+\log\bigg\{\frac{(2r)^{r/2-1/2}}{\sqrt{2}\pi^{(r+1)/2}}\int_{\mathcal {M}}\|D M_s\|_rds\bigg\}\bigg],
\end{align}
where $M_s$ is a $d\times r$ matrix with orthonormal columns spanning ${\cal T}_s\mathcal {M}.$ Then
\begin{align*}
\lim_{h\rightarrow0}\mathbb{P}\Big\{\sup_{t\in\mathcal {M}}|Z_h(t)|\leq \phi(z) \Big\}=\exp\{-2\exp\{-z\}\}.
\end{align*}
\end{theorem}

\vspace*{0.3cm}

This result will play a key role in the proof of Theorem~\ref{LevelSetConf}, which is presented now. First we are going to prove (\ref{ContourConf}) for $d\geq2$. Recall $\|K\|_2^2 = \int K^2(u)du$. Let 
\begin{align*}
Y_n(x)=\frac{\sqrt{nh^d}(\widehat f ( x ) -  f (x ) - \wh\beta(x))}{\sqrt{\|K\|_2^2\,f(x)}} = \frac{\sqrt{nh^d}(\widehat f ( x ) -  \mathbb{E}\widehat f (x ) + \beta(x) - \wh\beta(x))}{\sqrt{\|K\|_2^2\,f(x)}}.
\end{align*}
Let
\begin{align*}
b(z) = \sqrt{2(d-1)\log{h^{-1}}}&+\frac{1}{\sqrt{2(d-1)\log{h^{-1}}}}\bigg[z+\left(\frac{d}{2}-1\right)\log{\log{h^{-1}}}\nonumber\\
&\hspace{2cm}+\log\bigg\{\frac{(2d-2)^{d/2-1} s_K^{d-1}}{\sqrt{2}\pi^{d/2}}\mathscr{V}_{d-1}(\mathcal M)\bigg\}\bigg].
\end{align*}
To prove (\ref{ContourConf}), using Slutsky's Theorem, it suffices to show
\begin{align}\label{Ynextreme}
\lim_{n\rightarrow \infty}\mathbb{P}\left\{\sup_{x\in\mathcal{M}} | Y_n(x) | \leq b(z)\right\}=\exp\{-2\exp\{-z\}\},
\end{align}
because we have
\begin{align}\label{ciequivalence}
&\lim_{n\rightarrow \infty}\mathbb{P}\Big\{\sup_{x\in\mathcal{M}}\frac{\sqrt{nh^d}|\widehat f(x) -  \wh\beta(x) - c| }{\sqrt{\|K\|_2^2c}}\leq b(z)\Big\} \nonumber\\
=&\lim_{n\rightarrow \infty}\mathbb{P}\Bigg\{|\widehat f(x) -  \wh\beta(x) - c|\leq \frac{ b(z)\sqrt{\|K\|_2^2c}}{\sqrt{nh^d}}, \;\; \forall\; x\in\mathcal{M}\Bigg\} \nonumber\\
%
%
%
=&\lim_{n\rightarrow \infty}\mathbb{P}\Bigg\{\mathcal{M}\subset (\widehat f^{bc})^{-1}\left[c-\frac{ b(z)\sqrt{\|K\|_2^2c}}{\sqrt{nh^d}},c+\frac{ b(z)\sqrt{\|K\|_2^2c}}{\sqrt{nh^d}}\right]\Bigg\}.
\end{align}

To prove (\ref{Ynextreme}) we show the following two properties:
\begin{align}\label{stochasterm}
\lim_{n\rightarrow \infty}\mathbb{P}\left\{\sup_{x\in\mathcal{M}} \left | \frac{\sqrt{nh^d}(\widehat f ( x ) -  \mathbb{E}\widehat f (x ) }{\sqrt{\|K\|_2^2\,c}} \right | \leq b(z)\right\}=\exp\{-2\exp\{-z\}\}
\end{align}
and
\begin{align}\label{biasneg}
\sqrt{\log{h^{-1}}} \sup_{x\in\mathcal{M}} \left| \frac{\sqrt{nh^d}( \beta(x) - \wh\beta(x))}{\sqrt{\|K\|_2^2\,c}} \right| = o_p(1).
\end{align}

Using the uniform convergence rates for kernel density derivatives (see Lemma 3 in Arias-Castro et al. 2016 ) we obtain
\begin{align*}
&\textstyle{\sup\limits_{x\in\mathcal{M}} \left| \frac{\sqrt{nh^d}( \wh\beta(x) - \mathbb{E} \wh\beta(x) )}{\sqrt{\|K\|_2^2\,c}} \right| = O_p\left( \sqrt{nh^{d+4}} \beta_{n,l}^{(2),E} \right)},\quad\\
\text{and} \quad
%
&\textstyle{\sup\limits_{x\in\mathcal{M}} \left| \frac{\sqrt{nh^d}(  \mathbb{E} \wh\beta(x) - \beta(x) )}{\sqrt{\|K\|_2^2\,c}} \right| = O( \sqrt{nh^{d+4}}(l^2+h^2) ).}
\end{align*}
Property (\ref{biasneg}) now follows by using assumptions (\textbf{H1})$^0$ and (\textbf{H2})$^0$. Next we will show (\ref{stochasterm}). With
$\mathcal{F} = \big\{g_x(y)=\frac{1}{\sqrt{h^d\|K\|_2^2 c}}K\left(\frac{x-y}{h}\right):\;x\in\mathcal{M}\big\},$
and
$\mathbb{G}_n(g) = \frac{1}{\sqrt{n}}\sum_{i=1}^n\left(g(X_i)-\mathbb{E}g(X_1)\right),$  $\forall g\in\mathcal{F}.$
%
we can write 
%
$\frac{\sqrt{nh^d}(\widehat f(x)-\mathbb{E}\widehat f(x))}{\sqrt{\| K\|_2^2 c}} = \mathbb{G}_n(g_x),x\in\mathcal{M}$
and thus
\begin{align*}
\sup_{x\in\mathcal{M}} \left| \frac{\sqrt{nh^d}(\widehat f(x)-\mathbb{E}\widehat f(x))}{\sqrt{\| K\|_2^2 c}} \right| = \sup_{g_x\in\mathcal{F}} \left|\mathbb{G}_n(g_x)\right|.
\end{align*}

Let $\mathbb{B}$ be a centered Gaussian process on $\mathcal{F}$ such that for all $g_x, g_y\in\mathcal{F}$, 
$\mathbb{E}(\mathbb{B}(g_x)\mathbb{B}(g_y))=$  $ {\rm Cov}(g_x(X_1), g_y(X_1)).$ 
%
Applying Corollary 2.2 in Chernozhukov et al. (2014) we have that for all $\gamma\in(0,1)$ and $n$ sufficiently large 
\begin{align}\label{ChernozhukovApp}
&\textstyle{\mathbb{P}\left(\Big| \sup\limits_{x\in\mathcal{M}} \Big| \frac{\sqrt{nh^d}(\widehat f(x)-\mathbb{E}\widehat f(x))}{\sqrt{\| K\|_2^2 c}} \Big| - \sup\limits_{g\in\mathcal{F}}|\mathbb{B}(g)|\Big| \right.}\nonumber\\
&\hspace*{0.3cm}\textstyle{\left.>A_1\frac{\log^{2/3}(n)}{\gamma^{1/3}(nh^d)^{1/6}} + A_2\frac{\log^{3/4}(n)}{\gamma^{1/2}(nh^d)^{1/4}} + A_3\frac{\log(n)}{\gamma^{1/2}(nh^d)^{1/2}}\right)} \textstyle{\leq A_4\Big(\gamma + \frac{\log(n)}{n}\Big).}
\end{align}
where $A_1, A_2, A_3$ and $A_4$ are some constants. See Proposition 3.1 in Chernozhukov et al. (2014) for a similar derivation. 
%
%
%
%
Since $ \mathbb{E} \left[ \sup_{g\in\mathcal{F}}|\mathbb{B}(g)| \right] = O(\sqrt{\log{n}})$ (by Dudley's inequality for Gaussian processes, c.f. van der Vaart and Wellner, 1996, Corollary 2.2.8), with the choice of $\gamma = 1/\log{n}$, we apply Lemma 2.4 in Chernozhukov et al. (2014) and have
\begin{align*}
\textstyle{\sup\limits_{t} \left|\mathbb{P}\left(\sup_{x\in\mathcal{M}} \left| \frac{\sqrt{nh^d}(\widehat f(x)-\mathbb{E}\widehat f(x))}{\sqrt{\| K\|_2^2 c}} \right|<t\right) - \mathbb{P}\left(\sup\limits_{g\in\mathcal{F}}|\mathbb{B}(g)|<t\right)\right| = o(1).}
\end{align*}
It remains to show that
$\lim_{h\rightarrow 0} \mathbb{P}\left(\sup_{g\in\mathcal{F}}|\mathbb{B}(g)|< b(z) \right) = \exp\{-2\exp\{-z\}\}.$
%
Let $W$ and $B$ be $d$-dimensional Wiener process and Brownian bridge, respectively. Put 
$$U(x) = \frac{1}{h^{d/2} \| K\| \sqrt{f(x)}} \int_{\mathbb{R}^d} K \left( \frac{x-s}{h} \right) dB(M(s)),$$ where $M$ is the Rosenblatt transformation (cf. Rosenblatt, 1976).
Then
\begin{align*}
\sup_{g\in\mathcal{F}} \left|\mathbb{B}(g)\right| \stackrel{d}{=} \sup_{x\in\mathcal{M}} |U(x)|.
\end{align*}

Let further
$\widetilde U(x) = \frac{1}{h^{d/2}\| K \|} \int K\left( \frac{x-s}{h} \right) dW(s).$ 
%
Following the arguments on page 1013 of Rosenblatt (1976) (also see Proposition 2.2 in Bickel and Rosenblatt, 1973), we have 
$\sup_{x\in\mathcal{M}} |U(x) - \widetilde U(x)| = O_p(h^{1/2}).$ 
%
We then only need to show
\begin{align}\label{Uextreme}
\lim_{h\rightarrow 0} \mathbb{P}\left(\sup_{x\in\mathcal{M}} | \widetilde U(x)|< b(z) \right) = \exp\{-2\exp\{-z\}\}.
\end{align}

Next we are going to apply the probability results in Theorem \ref{ProbabilityResult}. Under our assumptions the isosurface $\mathcal{M}$ is a $d-1$ dimensional $C^1$ submanifold in $\mathbb{R}^d$ (see Theorem 2 in Walther, 1997). As discussed in Remark~\ref{remark}b), the reach of the isosurface is positive under our assumptions. It is easy to verify that the conditions for $Q(\delta)$ in (\ref{SupGauss2}) since $K$ is assumed to have bounded support. We will verify (\ref{varfunc}) and (\ref{SupGauss1}) in what follows.\\[5pt]
%
%
%
%
%
%
First observe that, as $\Delta x/h\rightarrow0$,
\begin{align}\label{covapprox}
Cov( \widetilde U(x), \widetilde U(x+\Delta x)) & = 1 - h^{-2} \Delta x^T\Sigma \Delta x + o(h^{-2} \|\Delta x\|^2),
\end{align}
where $\Sigma$ is a $d\times d$ symmetric matrix with the $(i,j)$-th element
\begin{align}\label{SigmaNotation}
\Sigma_{i,j} = \frac{\int \frac{\partial K(u)}{\partial u_i}\frac{\partial K(u)}{\partial u_j}du}{2\|K\|_2^2}.
\end{align}
Note that the little o term in (\ref{covapprox}) is uniform in $x\in f^{-1}[c-\delta_0,c+\delta_0]$ and $h\in(0,1]$, where $\delta_0$ appears in assumption (F2). Let $\|\partial K\|_2^2 =  \int \left(\frac{\partial K(u)}{\partial u_1}\right)^2du$. Then due to the symmetry of $K$, $\Sigma_{i,j} = s_K^2\delta_{i,j}$ where $s_K^2 = \frac{\|\partial K\|_2^2 }{2\|K\|_2^2 }$ and $\delta_{i,j}$ is the Kronecker delta. That is $\Sigma=s_K^2\bf{I}$. Then (\ref{varfunc}) and (\ref{SupGauss1}) are satisfied and we can apply Theorem \ref{ProbabilityResult} to get (\ref{Uextreme}) by showing that $b(z)$ is equivalent to
\begin{align}
&\textstyle{\sqrt{2(d-1)\log{\frac{1}{h}}}+\frac{1}{\sqrt{2(d-1)\log\{\frac{1}{h}}\}}\Big[z+\left(\frac{d}{2}-1\right)\log\log\frac{1}{h}}\nonumber\\
&\hspace*{4cm}
+\textstyle{\log\Big\{\frac{(2d-2)^{d/2-1}}{\sqrt{2}\pi^{(d/2)}}\int_{\mathcal {M}}\|\Sigma^{1/2} M_s\|_{d-1}ds\Big\}\Big],}\label{bzexpress}
\end{align}
where $M_s$ is a $d\times (d-1)$ matrix with orthonormal columns spanning ${\cal T}_s\mathcal {M}$, and $\|\cdot\|_{d-1}$ is the sum of squares of all minors of order $d-1$ for a $d\times (d-1)$ matrix. Note that since $\Sigma = s_K^2\bf{I}$ we have $\|\Sigma^{1/2} M_s\|_{d-1} = s_K^{(d-1)} \|M_s\|_{d-1}$. By the Cauchy-Binet formula (cf. Broida and Williamson 1989, pp 208-214), $\|M_s\|_{d-1} = \sqrt{det(M_s^TM_s)}=1$. Recall that $\mathscr{V}_{d-1}$ is the $(d-1)$-dimensional Hausdorff measure. The integral in (\ref{bzexpress}) can be simplified as 
\begin{align*}
\int_{\mathcal {M}}\|\Sigma^{1/2} M_s\|_{d-1}ds = s_K^{(d-1)} \mathscr{V}_{d-1} (\mathcal {M}).
\end{align*}

Now we have verified the expression in (\ref{bzexpress}) is equivalent to $b(z)$ and therefore (\ref{ContourConf}) follows.\\

Next we prove the case $d=1$. Following the discussion after the assumptions, denote $\mathcal{M} = \{x_i, i=1,\cdots, N\}$. Following similar argument at the beginning of the proof for $d\geq 2$ and using results from Bickel and Rosenblatt (1973), we have that the asymptotic distribution of 
%
$\textstyle{\sup\limits_{x\in\mathcal{M}}\frac{\sqrt{nh}|\widehat f^{bc}(x) - c |}{\sqrt{\|K\|_2^2c}}}$
%
is the same as that of 
%
$\sup_{x\in\mathcal{M}}|\tilde U(x)|,$
%
as $n\rightarrow\infty$. When $h$ is small enough, this supremum becomes $\max_{i=1,\cdots, N}Z_i$, where $Z_i$'s are i.i.d. standard normal random variables. Therefore, as $n\rightarrow\infty$,
\begin{align}
\textstyle{\sup\limits_{x\in\mathcal{M}}\frac{\sqrt{nh}|\widehat f^{bc}(x) - c |}{\sqrt{\|K\|_2^2c}} \overset{d}\longrightarrow \max_{i=1,\cdots, N}Z_i.}\label{AsymptoticFD1}
\end{align}
Recall that $\Phi$ is the standard normal c.d.f. Then the c.d.f of $\max_{i=1,\cdots, N}Z_i$ is $\Phi^N$. By Theorem 3.1 in Biau et al. (2007), $\wh N =N$ for $n$ large enough with probability 1. Following the same argument as for (\ref{ciequivalence}), we obtain (\ref{ContourConf}) for $d=1$. \hfill $\square$

\subsection{Proof of Corollary~\ref{UpperLSConf}}
The proof for $d=1$ is trivial. Now we show the proof for $d\geq2$. Since
%
\begin{align*}
&\mathbb{P}\left\{ \widehat{C}_{n,1}^+(1-\alpha) \subset \mathcal{L} \subset \widehat{C}_{n,1}^-(1-\alpha) \right\} \\
=\; & \mathbb{P}\left\{ \mathcal{M} \subset \widehat{C}_{n,1}(1-\alpha) \right\}
\mathbb{P}\left\{ \widehat{C}_{n,1}^+(1-\alpha) \subset \mathcal{L} \subset \widehat{C}_{n,1}^-(1-\alpha) \; | \;  \mathcal{M} \subset \widehat{C}_{n,1}(1-\alpha) \right\},
\end{align*}
using Theorem~\ref{LevelSetConf}, we only need to show that
\begin{align}\label{ConditionalProb}
\lim_{n\rightarrow\infty}\mathbb{P}\left\{ \widehat{C}_{n,1}^+(1-\alpha) \subset \mathcal{L} \subset \widehat{C}_{n,1}^-(1-\alpha) \; | \;  \mathcal{M} \subset \widehat{C}_{n,1}(1-\alpha) \right\} = 1.
\end{align}
%
In fact, we can obtain the following result: for all $L>0$ as $n\rightarrow\infty$,
\begin{align}\label{ConditionalProbStr}
\mathbb{P}\left\{ \widehat{C}_{n,1}^+(1-\alpha) \subset \mathcal{L} \subset \widehat{C}_{n,1}^-(1-\alpha) \; | \;  \mathcal{M} \subset \widehat{C}_{n,1}(1-\alpha) \right\} = 1 - O(n^{-L}).
\end{align}
Note that the event $\{ \widehat{C}_{n,1}^+(1-\alpha) \subset \mathcal{L} \subset \widehat{C}_{n,1}^-(1-\alpha) \}$ is equivalent to $E_1\cap E_2$, where 
$E_1 = \{ f(y)\geq c$, $\forall y$ s.t. $\widehat f^{\,bc}(y) \geq c+\wh a_{1-\alpha}^{(d)}$ $\}$, and $E_2 = \{$ $\widehat f^{\,bc}(y) \geq c - \wh a_{1-\alpha}^{(d)}$ $\forall y$ s.t. $f(y) \geq c$ $\}$.
We also denote $E_0 = \{ \mathcal{M} \subset \widehat{C}_{n,1}(1-\alpha) \} = \{$ $c -  \wh a_{1-\alpha}^{(d)} \leq \widehat f^{\,bc}(x) \leq c +  \wh a_{1-\alpha}^{(d)}$, $\forall x$ s.t. $f(x)=c$ $\}$. Then (\ref{ConditionalProbStr}) can be written in the form of $\mathbb{P}(E_1\cap E_2 | E_0) = 1- O(n^{-L})$ or equivalently, 
\begin{align}\label{transformedstate}
\mathbb{P}(E_1^\complement\cup E_2^\complement | E_0)  = O(n^{-L}),
\end{align}
where 
$E_1^\complement = \{$ $\exists y$ s.t. $f(y) < c \;\&$ $\widehat f^{\,bc}(y) \geq c+\wh a_{1-\alpha}^{(d)}$ $\}$, and 
$E_2^\complement = \{$ $\exists y$ s.t. $f(y) \geq c \;\&$ $\widehat f^{\,bc}(y) < c - \wh a_{1-\alpha}^{(d)}$  $\}$.
Property (\ref{transformedstate}) obviously follows from 
\begin{align}\label{ToShow1}
    \mathbb{P}(E_1^\complement | E_0) = O(n^{-L})
\end{align} 
and 
\begin{align*}
\mathbb{P}(E_2^\complement | E_0) = O(n^{-L}).
\end{align*} 
We only show (\ref{ToShow1}). To this end, we first introduce two more events. For some $C>0$ large enough, define
\begin{align*}
&\textstyle{B_0 = \Big\{\sup_{z\in\mathbb{R}^d} |\wh f^{bc}(z) - f(z)| \leq C\beta_{n,h}^{(0)} \Big\}},\quad\\
&\textstyle{B_1 = \Big\{\sup_{z\in\mathbb{R}^d} \|\nabla \wh f^{bc}(z) - \nabla f(z)\| \leq C\beta_{n,h}^{(1)} \Big\},}
\end{align*}
where we use the notation introduced in (\ref{beta-notation}). It follows from the proof on page 207 of Mammen and Polonik (2013) that $\mathbb{P}(B_0) = 1- O(n^{-L})$. Also note that
\begin{align*}
\sup_{z\in\mathbb{R}^d} \|\nabla \wh f^{\,bc}(z) - \nabla f(z)\| \leq \sup_{z\in\mathbb{R}^d} \|\nabla \wh f(z) - \mathbb{E} \nabla \wh f(z)\| + \sup_{z\in\mathbb{R}^d} \|\nabla \wh \beta(z) - \nabla \beta (z)\|.
\end{align*}
It is known (e.g. see Theorem 1 in Einmahl and Mason 2005, and Lemma 3 in Arias-Castro et al. 2016 ) that 
\begin{align}
& \sup_{z\in\mathbb{R}^d} \|\nabla \wh f(z) - \mathbb{E} \nabla \wh f(z)\| = O_{a.s.}\left( \beta_{n,h}^{(1),E}\right),\label{Einmahl1}\\
 &\sup_{z\in\mathbb{R}^d} \|\nabla \wh \beta(z) - \nabla \beta (z)\| = O_{a.s.}\left( h^2 \beta_{n,l}^{(3)} + h^2\right).\label{Einmahl2}
\end{align}
Following an argument similar to that in Mammen and Polonik (2013), we thus obtain $\mathbb{P}(B_1) = 1- O(n^{-L})$. These rates of convergences of $P(B_0)$ and $P(B_1)$ will be used in the following.\\

With $\delta_0$ given in Assumption (\textbf{F2}), let $\mathcal{M}^{\delta_0} =\{x: c-\delta_0\leq f(x) \leq c+\delta_0\}$. We now split up $E_1^\complement$ into 
\begin{align*}
E^\complement_{11} &= \{\exists y \in {\cal M}^{\delta_0} \text{s.t. }\, f(y) < c \;\& \widehat f^{\,bc}(y) \geq c+\wh a_{1-\alpha}^{(d)} \}\\
E^\complement_{12} &= \{\exists y \notin {\cal M}^{\delta_0} \text{s.t. }\, f(y) < c \;\& \widehat f^{\,bc}(y) \geq c+\wh a_{1-\alpha}^{(d)} \}.
\end{align*}
Note that these two sets are disjoint and $E_1^\complement = E^\complement_{11} \cup E^\complement_{12}.$  To show (\ref{ToShow1}), we now show that both $P(E_{11}^\complement| E_0) = O(n^{-L})$ and $P(E_{12}^\complement| E_0) = O(n^{-L}).$\\

First we show that $E_{12}^\complement \cap E_0 \cap B_0 = \emptyset$ for large enough $n$. To this end, let $x\in\mathcal{M}$, i.e. $f(x) = c$. Then, on $E_0$, we have $c -  \wh a_{1-\alpha}^{(d)} \leq \widehat f^{\,bc}(x) \leq c +  \wh a_{1-\alpha}^{(d)}$.  Let $y$ be a point that makes $E_{12}^\complement$ occur, i.e., $y\notin \mathcal{M}^{\delta_0}$ and $f(y) < c$ and $\widehat f^{\,bc}(y) \geq c+\wh a_{1-\alpha}^{(d)}$. 
%
Since $y\notin \mathcal{M}^{\delta_0},$ we in fact have $f(y)<c-\delta_0$. Note that $$  [\widehat f^{\,bc}(y) - \widehat f^{\,bc}(x)] - [f(y)-f(x)] \leq |\widehat f^{\,bc}(y) - f(y)| + |\widehat f^{\,bc}(x) - f(x)|,$$ 
which implies that, on $B_0$, 
\begin{align*}
\widehat f^{\,bc}(y) - \widehat f^{\,bc}(x) \leq |\widehat f^{\,bc}(y) - f(y)| + |\widehat f^{\,bc}(x) - c| +f(y) -c < -\delta_0 + 2C\beta_{n,h}^{(0)}.
\end{align*}
Also notice that, on $B_0,$ we have $\widehat f^{\,bc}(x) \leq c + C\beta_{n,h}^{(0)}$, and therefore $\widehat f^{\,bc}(y)< c-\delta_0 + 3C\beta_{n,h}^{(0)}$, which, for large enough $n$ cannot occur when $\widehat f^{\,bc}(y) \geq c+\wh a_{1-\alpha}^{(d)}$. Therefore we get $E_{12}^\complement \cap E_0 \cap B_0 = \emptyset$ for $n$ large enough.\\[5pt]
Next we show that $E_{11}^\complement \cap E_0 \cap B_1 = \emptyset$ for large enough $n$. So we now assume that $y\in\mathcal{M}^{\delta_0}$. Consider the integral curve $\X_y(t)$,  which is driven by $\nabla f/\|\nabla f\|^2$, starting from $y$. Recall that $\|\nabla f(z)\| > \epsilon_0>0$ for $z\in\mathcal{M}^{\delta_0}$ by assumption {\bf (F2)}. Let $\theta = c - f(y) > 0$ (on $E_{12}^\complement$). Using the property of the integral curve described in (\ref{gradientpathproperty}), we have $\X_y(\theta) \in \mathcal{M}$.  
%
On $B_1,$ we have $\langle \nabla \wh f^{\,bc}(z), \nabla f(z)/\|\nabla f(z)\|^2 \rangle>0,\; \forall z\in \{\X_y(t): \; t\in[0,\theta]\}$ for $n$ large enough. This means $\widehat f^{\,bc}$ keeps increasing on the trajectory of $\{\X_y(t): \; t\in[0,\theta]\}$. Therefore, $\widehat f^{\,bc}(\X_y(\theta))> \widehat f^{\,bc}(y) \geq c+\wh a_{1-\alpha}^{(d)}$, which contradicts $E_0$. Therefore, $E_{11}^\complement \cap E_0 \cap B_1 = \emptyset$ for $n$ large enough.\\[3pt] 
This now results in the following. For $n$ large enough
\begin{align*}
\mathbb{P}(E_1^\complement| E_0)&
=[\mathbb{P}(E_{11}^\complement \cap E_0)+ \mathbb{P}(E_{12}^\complement \cap E_0)]/ \mathbb{P}(E_0) \\
=\; & \big[\mathbb{P}(E_{11}^\complement \cap E_0 \cap B^\complement_1) + \mathbb{P}(E_{12}^\complement  \cap E_0 \cap B_0^\complement)\big] / \mathbb{P}(E_0) \\
\le\; & [\mathbb{P}(B_1^\complement) + \mathbb{P}(B_0^\complement)]/\mathbb{P}(E_0).
\end{align*}
Since both $\mathbb{P}(B_1^\complement) = O(n^{-L})$ and $\mathbb{P}(B_0^\complement) = O(n^{-L}),$ and $P(E_0) \to 1-\alpha$, the assertion follows. The fact that $\mathbb{P}(E_2^\complement| E_0) = O(n^{-L})$ can be shown in a similar way. This completes the proof.
\hfill$\square$

\subsection{Proof of Theorem~\ref{smooth-bootstrap-method}}

Assuming (\ref{confme2}) (or (\ref{confm2})) is true, then (\ref{confle2}) (or (\ref{confl2})) can be proved in a very similar way as for Corollary~\ref{UpperLSConf}. In particular, notice that the key result (\ref{ConditionalProbStr}) in the proof of Corollary~\ref{UpperLSConf} can be replaced by (say, for the proof of (\ref{confl2}))
\begin{align*}
\mathbb{P}\left\{ \widehat{C}_{n,2}^{*,+}(1-\alpha) \subset \mathcal{L} \subset \widehat{C}_{n,2}^{*,-}(1-\alpha) \; | \;  \mathcal{M} \subset \widehat{C}_{n,2}(1-\alpha) \right\} = O(n^{-L}).
\end{align*}

Next we prove (\ref{confme2}) and (\ref{confm2}). We will use $C$ to denote a generic constant that may be different at different occurrences. 
For any $\widehat c>0$, the following three events are equivalent:
\begin{align*}
\mathcal{M}^E \subset \wh f^{-1}[c-\widehat c,\, c+\widehat c] \; &\Leftrightarrow \; \wh f(x) \in [c-\widehat c,\, c+\widehat c], \; \forall x\in\mathcal{M}^E \; \\
&\Leftrightarrow \; \sup_{x\in\mathcal{M}^E} |\wh f(x) -\mathbb{E}\wh f(x)|\leq \widehat c.
\end{align*}
Similarly
\begin{align*}
\mathcal{M} \subset \wh f^{-1}[c-\widehat c,\, c+\widehat c] \; &\Leftrightarrow \; \wh f(x) \in [c-\widehat c,\, c+\widehat c], \; \forall x\in\mathcal{M} \; \\
&\Leftrightarrow \; \sup_{x\in\mathcal{M}} |\wh f(x) - f(x)|\leq \widehat c.
\end{align*}
Therefore it suffices to show that
\begin{align}\label{supexpdiffrate}
\textstyle{\mathbb{P}\Big( \sup\limits_{x\in\mathcal{M}^E} |\wh f(x) -\mathbb{E}\wh f(x)|\leq \widehat c_{1-\alpha}^{*,E} \Big) = (1-\alpha) + O\big(\Psi_n(\gamma_n^E)\big),}
\end{align}
and\\[-25pt]
\begin{align}\label{supdiffrate}
\textstyle{\mathbb{P}\Big( \sup\limits_{x\in\mathcal{M}} |\wh f(x) - f(x)|\leq \widehat c_{1-\alpha}^{*} \Big) = (1-\alpha) + \big(\Psi_n(\gamma_n)\big).}
\end{align}

The proofs for these two results are similar. We first show (\ref{supexpdiffrate}) and then briefly sketch the proof for (\ref{supdiffrate}). It is known from page 209 in Mammen and Polonik (2013) (also see Theorem 3.1 in Neumann (1998)) that for some $C<\infty$
\begin{align}\label{expcomparison}
&\mathbb{P} \Big( \sup_{x\in\mathbb{R}^d} |\wh f(x) -\mathbb{E}\wh f(x) - (\widehat f^*(x) - \mathbb{E}^* \widehat f^*(x))| > C \Big(\beta_{n,h}^{(0),E} +\sqrt{\beta_{n,g}^{(0)}} \,\Big) \beta_{n,h}^{(0),E} \Big)\nonumber \\
&\hspace*{2cm}= O(n^{-L}),
\end{align}
for an arbitrarily large $L>0$, which implies
\begin{align}\label{supexpdiff}
&\textstyle{\mathbb{P} \Big( \Big|\sup\limits_{x\in\wh{\mathcal{M}}} |\wh f(x) -\mathbb{E}\wh f(x)| - \sup\limits_{x\in\wh{\mathcal{M}}} |\widehat f^*(x) - \mathbb{E}^* \widehat f^*(x)| \Big| > C \Big(\beta_{n,h}^{(0),E} +\sqrt{\beta_{n,g}^{(0)}}\,\Big)\Big)}\nonumber \\
&\hspace*{2cm}= \;O(n^{-L}).
\end{align}

When $n$ is large enough there exists $C_1>0$ such that $\wh{\mathcal{M}} \subset \bigcup_{x\in\mathcal{M}^E} \mathscr{B}(x, C_1\beta_n^E)$, where $\mathscr{B}(x,r)$ is the ball in $\mathbb{R}^d$ with center $x$ and radius $r$. Therefore
\begin{align*}
&\Big| \sup_{x\in\wh{\mathcal{M}}} |\wh f(x) -\mathbb{E}\wh f(x)| - \sup_{x\in\mathcal{M}^E} |\wh f(x) -\mathbb{E}\wh f(x)|\Big|\\
 \leq &\sup_{\|x-y\|\leq C_1\beta_{n,h}^{(0),E}} \left| (\wh f(x) -\mathbb{E}\wh f(x)) - (\wh f(y) -\mathbb{E}\wh f(y))\right|.
\end{align*}

It follows from an argument similar to the one given in Mammen and Polonik (2013), page 209,  that, for some $C>0$,
\begin{align*}
\textstyle{\mathbb{P}\Big( \sup\limits_{\|x-y\|\leq C_1\beta_{n,h}^{(0),E}} \left| (\wh f(x) -\mathbb{E}\wh f(x)) - (\wh f(y) -\mathbb{E}\wh f(y))\right| \geq C \beta_{n,h}^{(1),E}\,\beta_{n,h}^{(0),E} \Big) = O(n^{-L}).}
\end{align*}
This then leads to 
\begin{align}\label{Mexpdiff}
\textstyle{\mathbb{P}\Big( \Big| \sup\limits_{x\in\wh{\mathcal{M}}} |\wh f(x) -\mathbb{E}\wh f(x)| - \sup\limits_{x\in\mathcal{M}^E} |\wh f(x) -\mathbb{E}\wh f(x)|\Big| \geq C \beta_{n,h}^{(1),E}\beta^{(0),E}_{n,h} \Big) = O(n^{-L}),}
\end{align}
which combining with (\ref{supexpdiff}) further implies
\begin{align}\label{lemma24A1}
\textstyle{\mathbb{P}\Big( \Big| \sup\limits_{x\in\wh{\mathcal{M}}} |\widehat f^*(x) - \mathbb{E}^* \widehat f^*(x)| - \sup\limits_{x\in\mathcal{M}^E} |\wh f(x) -\mathbb{E}\wh f(x)|\Big| \geq C\gamma_n^E  \Big) = O(n^{-L}),}
\end{align}
where $\gamma_n^E$ is given in (\ref{gammanE}). As a result of Proposition 3.1 in Neumann (1998) we obtain with $\Psi_n$ as in (\ref{Neumann}),
\begin{align}\label{neumannres}
\mathbb{P} \Big( \sup_{x\in\mathcal{M}^E} |\wh f(x) -\mathbb{E}\wh f(x)| \in [c, d]\Big) = O\big( \Psi_n(d-c)\big),
\end{align}
where this rate holds uniformly in $0 \le c < d < \infty$ Thus, for some $C>0$, we have
\begin{align}\label{lemma24A2}
\sup_{t\in\mathbb{R}} \mathbb{P} \Big( \sup_{x\in\mathcal{M}^E} |\wh f(x) -\mathbb{E}\wh f(x)| \in [t, t+\gamma_n^E]\Big) \leq C\Psi_n\big(\gamma_n^E\big).
\end{align}
With (\ref{lemma24A1}) and (\ref{lemma24A2}), then (\ref{supexpdiffrate}) follows from Lemma 2.4 in Mammen and Polonik (2013).\\[4pt]
%
Next, we briefly outline the proof of (\ref{supdiffrate}). Following the proof on page 207 and Mammen and Polonik (2013) (also see Lemma 3 in Arias-Castro et al., 2016), we have that under our assumption, for some $C>0$,
\begin{align}\label{2ndderivrate}
\mathbb{P}\Big( \sup_{x\in\mathbb{R}^d} \Big| \frac{\partial^2}{\partial x_i \partial x_j} \widehat f_g(x) - \frac{\partial^2}{\partial x_i \partial x_j}  f(x)\Big| \geq C\beta_{n,g}^{(2)}\Big) = O(n^{-L}),\; \text{for all } i,j=1,\cdots,d.
\end{align}
Mammen and Polonik (2013, page 210) show that, for all $\delta>0$, there exists $C>0$ such that
\begin{align}\label{expdiffcomp}
&\textstyle{\sup\limits_{\|x\| \leq \delta} |[\mathbb{E}^* \widehat f^*(x) - \mathbb{E}\widehat f(x)] - [\widehat f_g(x) - f(x)]|} \nonumber\\
&\hspace{4cm} \textstyle{\leq Ch^2 \sum\limits_{i,j=1}^d \sup\limits_{\|x\|\leq \delta+\sqrt{d}h}  \left| \frac{\partial^2}{\partial x_i \partial x_j} \widehat f_g(x) - \frac{\partial^2}{\partial x_i \partial x_j}  f(x)\right|\; \text{\rm a.s.}}
\end{align}

%
%
Since $\wh{\mathcal{M}} \subset \mathscr{B}(0,\delta)$ for large $\delta$ when $n$ is large enough with probability one, by (\ref{expcomparison}), (\ref{2ndderivrate})  and (\ref{expdiffcomp}) we have
%
\begin{align}\label{bootmhatdiff}
&\mathbb{P} \Big( \Big|\sup_{x\in\wh{\mathcal{M}}} |\wh f(x) -  f(x)| - \sup_{x\in\wh{\mathcal{M}}} |\widehat f^*(x) -  \widehat f_g(x)| \Big|  \nonumber\\
&\hspace{2cm} > C \big( \big(\beta_{n,h}^{(0),E}+\beta_{n,g}^{(0)}\big) \beta_{n,h}^{(0),E} + h^2\beta_{n,g}^{(2)}\big)\Big) = O(n^{-L}).
\end{align}
Since there exists $C_2>0$ such that $\wh{\mathcal{M}} \subset \bigcup_{x\in\mathcal{M}^E} \mathscr{B}(x, C_2\beta^{(0)}_{n,h})$ when $n$ is large enough,  the following result similar to (\ref{Mexpdiff}) can be derived
\begin{align}\label{mmhatdiff}
\mathbb{P}\Big( \Big| \sup_{x\in\wh{\mathcal{M}}} |\wh f(x) - f(x)| - \sup_{x\in\mathcal{M}} |\wh f(x) - f(x)|\Big| \geq C \beta_{n,h}^{(1)} \beta^{(0)}_{n,h} \Big) = O(n^{-L}).
\end{align}
Then combining (\ref{bootmhatdiff}) and (\ref{mmhatdiff}) we get
\begin{align}\label{lemma24B1}
\mathbb{P}\Big( \Big| \sup_{x\in\wh{\mathcal{M}}} |\widehat f^*(x) - \widehat f_g(x)| - \sup_{x\in\mathcal{M}} |\wh f(x) - f(x)|\Big| \geq C\gamma_n  \Big) = O(n^{-L}),
\end{align}
where $\gamma_n$ is given in (\ref{gamman}). Mammen and Polonik (2013) modify (\ref{neumannres}) to 
\begin{align}\label{neumannres2}
\mathbb{P} \Big( \sup_{x\in\mathcal{M}} |\wh f(x) - f(x)| \in [c, d]\Big) = O\big( \Psi_n(d-c)\big),
\end{align}
which immediately gives, for some $C > 0$,
\begin{align}\label{lemma24B2}
\mathbb{P} \Big( \sup_{x\in\mathcal{M}} |\wh f(x) -  f(x)| \in [t, t+\gamma_n]\Big) \leq C\Psi_n\big(\gamma_n\big),
\end{align}
where (\ref{neumannres2}) and (\ref{lemma24B2}) hold uniformly over $0 \le c < d < \infty$ and $t \in \R$, respectively. Applying Lemma~\ref{MP} with $Z_n = \sup_{x\in\mathcal{M}} |\wh f(x) - f(x)|$ and $Z_n^* = \sup_{x\in\wh{\mathcal{M}}} |\widehat f^*(x) - \widehat f_g(x)|$, by using  (\ref{lemma24B1}) and (\ref{lemma24B2}), we obtain (\ref{supdiffrate}). \hfill $\square$

\subsection{Proof of Theorem~\ref{conf-region-explicit-bias}}
Following the same argument as in the proof of Theorem~\ref{smooth-bootstrap-method}, (\ref{cn3confl}) follows easily (using the proof of Corollary~\ref{UpperLSConf}) once (\ref{cn3confm}) is proved. We will only show the latter. Since\\[-20pt]
\begin{align*}
\beta(x) = \frac{1}{2} h^2\int u_1^2K(u)du \sum_{j=1}^d \frac{\partial^2}{\partial x_j \partial x_j}f(x) + O(h^4),
\end{align*}
\vspace*{-1cm}

we have
\begin{align}\label{betahatdiff}
\sup_{x\in\mathcal{M}} |\wh \beta(x) - \beta(x)| = O_{a.s.}(h^2\,\beta_{n,l}^{(2)}+h^4).
\end{align}
Similar to (\ref{2ndderivrate}), there exists $C>0$ such that, for an arbitrary $L>0$,
\begin{align}\label{2ndderivrate2}
\mathbb{P} \Big( \sup_{x\in\mathcal{M}} |\wh \beta(x) - \beta(x)| \geq C h^2\big(\beta_{n,l}^{(2)} +h^2 \big)\Big) = O(n^{-L}),
\end{align}
and, similar to (\ref{Mexpdiff}), we obtain
\begin{align}\label{Mexpdiff2}
\mathbb{P} \Big( \Big| \sup_{x\in\wh{\mathcal{M}}} |\wh f(x) -\mathbb{E}\wh f(x)| - \sup_{x\in\mathcal{M}} |\wh f(x) -\mathbb{E}\wh f(x)|\Big| \geq C \beta_{n,h}^{(1),E}\,\beta_{n,h}^{(0)} \Big) = O(n^{-L}).
\end{align}
Since 
%
 $\Big| \sup_{x\in\wh{\mathcal{M}}} |\wh f(x) -\mathbb{E}\wh f(x)| - \sup_{x\in\mathcal{M}} |\wh f(x) -\mathbb{E}\wh f(x)| \Big| + \sup_{x\in\mathcal{M}} |\wh \beta(x) - \beta(x)|\\
 %
  \geq  \Big| \sup_{x\in\wh{\mathcal{M}}} |\wh f(x) -\mathbb{E}\wh f(x)| -  \sup_{x\in\mathcal{M}} |\wh f^{\,bc}(x) - f(x) | \Big|,$
%
combining (\ref{2ndderivrate2}) and (\ref{Mexpdiff2}), we have
\begin{align}\label{Mexpdiff2res}
&\mathbb{P}\Big( \Big| \sup_{x\in\wh{\mathcal{M}}} |\wh f(x) -\mathbb{E}\wh f(x)| -  \sup_{x\in\mathcal{M}} |\wh f^{\,bc}(x) - f(x) | \Big| \geq C \Big( \beta_{n,h}^{(1),E}\,\beta_n + h^2(\beta_{n,l}^{(2)} +h^2) \Big) \Big) \nonumber\\[7pt]
&\hspace*{2cm}= O(n^{-L}).
\end{align}
By (\ref{supexpdiff}), we get
\begin{align}\label{lemma24C1}
\mathbb{P}\Big( \Big| \sup_{x\in\wh{\mathcal{M}}} |\wh f^*(x) -\mathbb{E}^*\wh f^*(x)| -  \sup_{x\in\mathcal{M}} |\wh f^{\,bc}(x) - f(x) | \Big| \geq C \gamma_n^{bc}\Big) = O(n^{-L}).
\end{align}
Similar to (\ref{neumannres2}), we have 
\begin{align}\label{neumannres3}
\mathbb{P} \Big( \sup_{x\in\mathcal{M}} |\wh f(x) - \mathbb{E} \wh f(x) | \in [c, d]\Big) = O\big( \Psi_n(d-c)\big)
\end{align}
where this convergence is uniform in $0 \le c < d < \infty$. Combining (\ref{2ndderivrate2}) and (\ref{neumannres3}), we have with $\kappa_n = Ch^2(\beta_{n,l}^{(2)}+h^2)$,
\begin{align}\label{neumannres4}
&\mathbb{P} \Big( \sup_{x\in\mathcal{M}} |\wh f^{\,bc}(x) - f(x) | \in [c, d]\Big) \nonumber\\
\leq &\, \mathbb{P} \Big( \sup_{x\in\mathcal{M}} |\wh f(x) - \mathbb{E} \wh f(x) | \in [c-\kappa_n, d+\kappa_n]\Big) + O(n^{-L}) \nonumber\\
= &\, O\Big( \Psi_n(d-c) + \kappa_n\sqrt{nh^d\log{n}}\Big).
\end{align}
This then implies that, for some $C>0$,
\begin{align}\label{lemma24C2}
\sup_{t\in\mathbb{R}} \mathbb{P} \Big( \sup_{x\in\mathcal{M}} |\wh f^{\,bc}(x) - f(x) | \in [t, t+\gamma_n^{bc}] \Big) \leq C\Psi_n(\gamma_n^{bc}).
\end{align}

Therefore, using Lemma~\ref{MP} with (\ref{lemma24C1}) and (\ref{lemma24C2}), we have
\begin{align*}
\mathbb{P}\Big( \sup_{x\in\mathcal{M}} |\wh f^{\,bc}(x) -f(x) |\leq \widehat c_{1-\alpha}^{*,E} \Big) = (1-\alpha) + O\big(\Psi_n(\gamma_n^{bc})\big),
\end{align*}
and the conclusion of the theorem follows. \hfill $\square$

\subsection{Proofs for Section 3}

We begin with two lemmas that will be needed in the proofs. Recall that $\wh f^{\,bc}(x) = \wh f (x) - \wh \beta(x)$ be the de-biased density estimator, and $\wh \X^{bc}_x$ be the integral curve driven by $\frac{\nabla \wh f^{\,bc}}{\|\nabla \wh f^{\,bc}\|^2}$ and $\wh\theta^{\,bc}_x$ be the corresponding time when $\wh\X^{bc}_x$ hits $\wh{\mathcal{M}}^{bc}$.

\begin{lemma}\label{existenceuniqueness}
Under assumptions (\textbf{F1}), (\textbf{F2}), (\textbf{K}), $\beta_{n,h}^{(1)}=o(1)$ and $\beta_{n,l}^{(3)}=o(1)$ as $n\rightarrow\infty$, with probability one we have that for $x\in\mathcal{M}$ the solution $\wh{\theta}_x^{bc}$ in (\ref{thetatilde}) exists and is unique for $n$ large enough. In such a case the mapping $x\mapsto \wh\X^{\,bc}_x(\wh\theta^{\,bc}_x)$ is bijective between $\mathcal{M}$ and $\wh{\mathcal{M}}^{bc}$.
\end{lemma}

{\bf Proof of Lemma \ref{existenceuniqueness}}
Throughout this proof we assume that the sample size is large enough. We first show that the solution $\wh{\theta}_x^{bc}$ in (\ref{thetatilde}) exists and is unique. Due to the strong consistency of $\wh f^{bc}$, we have that $\wh{\mathcal{M}}^{bc}\subset \mathcal{M}^{\delta_0} :=f^{-1}[c-\delta_0,c+\delta_0]$. Also due to the strong consistency of gradient estimator $\nabla\wh f^{bc}$ implied by (\ref{Einmahl1}) and (\ref{Einmahl2}), we have $\|\nabla\wh f^{bc}\|>\epsilon_0/2$ for all $x\in \mathcal{M}^{\delta_0}$ by assumption (\textbf{F2}). Since the trajectory of $\wh\X_x^{bc}$ is driven by $\nabla\wh f^{bc}$, for $x\in\mathcal{M}$ we have the existence and uniqueness of $\wh{\theta}_x^{bc}$ and $-\infty<\wh{\theta}_x^{bc}<\infty$.  Also using $\|\nabla\wh f^{bc}\|>0$ again, as a property of integral curves we have $\wh\X^{\,bc}_x(\wh{\theta}_{x_1}^{bc}) \neq \wh\X^{\,bc}_x(\wh{\theta}_{x_2}^{bc})$ if $x_1\neq x_2$ for $x_1,x_2\in\mathcal{M}$.\\[4pt]
The above argument also implies that the mapping $x\mapsto \wh\X^{\,bc}_x(\wh\theta^{\,bc}_x)$ is bijective between $\mathcal{M}$ and $\wh{\mathcal{M}}^{bc}$. The injectivity is an immediate consequence. Next we show the surjectivity. For any $y\in\wh{\mathcal{M}}^{bc}$, without loss of generality we assume $f(y)<c$. Since $\langle \nabla \wh f^{bc}, \nabla f/\|\nabla f\|^2\rangle >0 $ on $\mathcal{M}^{\delta_0}$, the value of $f$ keeps increasing on the trajectory of $\wh\X^{\,bc}_y$. Therefore there exists a finite time $\tilde\theta_y$ when $\wh\X^{\,bc}_y$ hits $\mathcal{M}$. Let $x=\wh\X^{\,bc}_y(\tilde\theta_y)$ and $\wh\theta^{\,bc}_x = - \tilde\theta_y$. Then $y=\wh\X^{\,bc}_x(\wh\theta^{\,bc}_x)$. \hfill$\square$
\begin{lemma}\label{EstimatedPath}
Under assumptions (\textbf{F1}), (\textbf{F2}), (\textbf{K}),  ({\textbf H1})$^2$ and ({\textbf H2})$^2$ we have with\\[-22pt]
\begin{align}
\tau_n = \beta_{n,h}^{(0),E} + h^2\left( \beta_{n,l}^{(2)}  + h^2\right),
\end{align}
%
%
that\\[-22pt]
\begin{align}\label{EstimatedPathTheta}
\sup_{x\in\mathcal{M}}|\wh\theta^{\,bc}_x| = \sup_{x\in\mathcal{M}}|\wh f^{\,bc}(x) - f(x)| =  O_{a.s.}\left( \tau_n \right),
\end{align}
\vspace*{-0.3cm}

\noindent
and
\vspace*{-0.3cm}

\begin{align}\label{EstimatedPathX}
\sup_{x\in\mathcal{M}}\left| \|\nabla\wh f^{\,bc}(\wh \X^{bc}_x(\wh \theta^{\,bc}_x))\|\|\wh \X^{\,bc}_x(\wh\theta^{\,bc}_x) -x\| - |\wh f^{\,bc}(x) - f(x)| \right| =O_{a.s.}\left( \tau_n^2 \right).
\end{align}
\end{lemma}
%

{\bf Proof of Lemma \ref{EstimatedPath}} 
By Lemma~\ref{existenceuniqueness}, for $x\in\mathcal{M}$, $\wh{\theta}_x^{bc}$ exists and is unique for large sample. For $x\in\mathcal{M}$, $ f(x)=c=\wh f^{\,bc}(\wh\X^{\,bc}_x(\wh\theta^{\,bc}_x))$ and therefore
\begin{align}\label{hatthetaX}
\wh\theta^{\,bc}_x = \wh f^{\,bc}(\wh\X^{\,bc}_x(\wh\theta^{\,bc}_x)) - \wh f^{\,bc}(x) = f(x) - \wh f^{\,bc}(x).
\end{align}
Consequently,
\begin{align*}
\sup_{x\in\mathcal{M}}|\wh \theta^{\,bc}_x| = \sup_{x\in\mathcal{M}}|\wh f^{\,bc}(x) - f(x)| = \sup_{x\in\mathcal{M}}|\widehat f(x) - \mathbb{E}\widehat f(x) + \beta(x) - \widehat \beta(x)| = O_{a.s.}\left( \tau_n \right),
%
\end{align*}
where we have used (\ref{betahatdiff}). This is (\ref{EstimatedPathTheta}). Next we prove (\ref{EstimatedPathX}). Without loss of generality we assume $\wh \theta^{\,bc}_x>0$ in what follows. We can write
\begin{align}\label{integralcurveapprox}
\wh \X^{\,bc}_x(\wh\theta^{\,bc}_x) -x &= \int_0^{\wh\theta^{\,bc}_x} \frac{\nabla\wh f^{\,bc}(\wh\X^{\,bc}_x(t))}{\|\nabla\wh f^{\,bc}(\wh\X^{\,bc}_x(t))\|^2}dt \nonumber\\
&= \frac{\nabla\wh f^{\,bc}(\wh\X^{\,bc}_x(\wh\theta^{\,bc}_x))}{\|\nabla\wh f^{\,bc}(\wh\X^{\,bc}_x(\wh\theta^{\,bc}_x))\|^2} \wh\theta^{\,bc}_x + \int_0^{\wh\theta^{\,bc}_x} \wh\eta^{\,bc}(t) dt.
\end{align}
where we denote $\wh\eta^{\,bc}(t) = \frac{\nabla\wh f^{\,bc}(\wh\X^{\,bc}_x(t))}{\|\nabla\wh f^{\,bc}(\wh\X^{\,bc}_x(t))\|^2} - \frac{\nabla\wh f^{\,bc}(\wh\X^{\,bc}_x(\wh\theta^{\,bc}_x))}{\|\nabla\wh f^{\,bc}(\wh\X^{\,bc}_x(\wh\theta^{\,bc}_x))\|^2}$. From (\ref{integralcurveapprox}) we obtain
\begin{align*}
\textstyle{\Big|  \|\wh \X^{\,bc}_x(\wh\theta^{\,bc}_x) -x\| - \frac{|\wh\theta^{\,bc}_x|}{\|\nabla\wh f^{\,bc}(\wh\X^{\,bc}_x(\wh\theta^{\,bc}_x))\|}    \Big| \leq  \sup\limits_{t\in[0,\wh \theta^{\,bc}_x]} \|\wh\eta^{\,bc}(t) \| |\wh\theta^{\,bc}_x|.}
\end{align*}
By (\ref{hatthetaX}) we have
\begin{align}\label{integralcurveapprox0}
& \textstyle{\sup\limits_{x\in\mathcal{M}}\big| \|\nabla\wh f^{\,bc}(\wh\X^{\,bc}_x(\wh\theta^{\,bc}_x))\|\|\wh \X^{\,bc}_x(\wh\theta^{\,bc}_x) -x\| - |f(x) - \wh f^{\,bc}(x)| \big|} \nonumber \\
& \hspace{1cm} \textstyle{\leq  \sup\limits_{x\in\mathcal{M}}\sup\limits_{t\in[0,\wh \theta^{\,bc}_x]} \|\wh\eta^{\,bc}(t) \| \sup_{x\in\mathcal{M}}|\wh\theta^{\,bc}_x| \sup_{x\in\mathcal{M}}\|\nabla\wh f^{\,bc}(\wh\X^{\,bc}_x(\wh\theta^{\,bc}_x))\| .}
\end{align}

Using assumptions ({\textbf H1})$^2$ and ({\textbf H2})$^2,$ we obtain for some small $\epsilon>0$
\begin{align*}
\sup_{y\in\mathcal{M}\oplus\epsilon} \max \{\| \nabla^2 \wh f^{\,bc}(y) -  \nabla^2 f(y)\|_F, \;\; \| \nabla \wh f^{\,bc}(y) - \nabla f(y)\| \}= o_{a.s.}(1),
\end{align*}
where $\|\cdot\|_F$ is the Frobenius norm and $\mathcal{M}\oplus\epsilon = \{x\in\mathbb{R}^d: \;d(x,\mathcal{M})\leq \epsilon\}$. 
By using Taylor expansion, we then have for some small $\epsilon>0$, $$\sup_{t\in[0,\wh \theta^{\,bc}_x]}\|\wh\eta^{\,bc}(t) \| \leq  \sup_{y\in\mathcal{M}\oplus\epsilon} \left[ \| \nabla^2 \wh f^{\,bc}(y)\|_F \| \nabla \wh f^{\,bc}(y)\|^{-1} \right] \wh \theta^{\,bc}_x.$$

Therefore using (\ref{integralcurveapprox0}) for some $C>0$, with probability one for $n$ large enough,
\begin{align*}
\sup_{x\in\mathcal{M}}\left| \|\nabla\wh f^{\,bc}(\wh\X^{\,bc}_x(\wh\theta^{\,bc}_x))\|\|\wh \X^{\,bc}_x(\wh\theta^{\,bc}_x) -x\| - |f(x) - \wh f^{\,bc}(x)| \right| \leq C \sup_{x\in\mathcal{M}} |\wh\theta^{\,bc}_x|^2.
\end{align*}
We conclude the proof of (\ref{EstimatedPathX}) by using (\ref{EstimatedPathTheta}). \hfill $\square$\\
%
%

{\bf Proof of Theorem \ref{asympt-Cn3-Cn4}.} 
First notice that (\ref{cn4confl}) (or (\ref{cn4starconfl})) follows from (\ref{cn4confm} (or (\ref{cn4starconfm})) using the same approach as in the proof of Corollary~\ref{UpperLSConf}. We point out that the events $E_1$ and $E_2$ defined in the proof of Corollary~\ref{UpperLSConf} need to be replaced by (say, for the proof of (\ref{cn4confl}))
\begin{align*}
&E_1 = \{f(y)\geq c, \forall y \text{ s.t. } y\in\wh C_{n,4}^+(1-\alpha)\},\\
\text{ and } &E_2 = \{ y\in\wh C_{n,4}^-(1-\alpha), \forall y \text{ s.t. } f(y)\geq c\}.
\end{align*}
Then we can show that 
\begin{align*}
&\mathbb{P}\left\{ \wh C_{n,4}^+(1-\alpha) \subset\mathcal{L}\subset\wh C_{n,4}^+ (1-\alpha) \; | \; \mathcal{M} \subset \wh C_{n,4}(1-\alpha) \right\} \\
= &\mathbb{P}\left\{ E_1\cap E_2 \;|\;  \mathcal{M} \subset \wh C_{n,4}(1-\alpha)\right\} = O(n^{-L}).
\end{align*}
Details are omitted. We focus on the proofs of (\ref{cn4starconfl}) and (\ref{cn4starconfm}) in what follows.\\[5pt]
{\bf Part 1.} The result immediately follows from Lemma \ref{EstimatedPath} and Theorem~\ref{LevelSetConf}.\\[5pt]
{\bf Part 2.} 
By (\ref{Mexpdiff2res}) and Lemma \ref{EstimatedPath} we obtain that for some $C>0$
\begin{align}\label{Mexpdiff2res2}
&\textstyle{\mathbb{P}\Big( \Big| \sup\limits_{x\in\wh{\mathcal{M}}} |\wh f(x) -\mathbb{E}\wh f(x)| -  \sup\limits_{x\in\mathcal{M}} \|\nabla\wh f^{\,bc}(\wh\X^{\,bc}_x(\wh\theta^{\,bc}_x))\|\|\wh \X^{\,bc}_x(\wh\theta^{\,bc}_x) -x\|\Big|} \nonumber\\
&\hspace{3cm}\textstyle{\geq C \big( \tau_n^2+\beta_{n,h}^{(1),E}\,\beta_{n,h}^{(0)} + h^2(\beta_{n,l}^{(2)}+h^2) \big) \Big) = O(n^{-L}).}
\end{align}
Noting that $\tau_n^2  =o\left( \beta_{n,h}^{(1),E}\,\beta_{n,h}^{(0)} + h^2\big(\beta_{n,l}^{(2)}+h^2\,\big) \right)$, the assertion follows from an argument similar to that in the proof of Theorem~\ref{conf-region-explicit-bias}. \hfill$\square$\\

{\bf The proof of Theorem \ref{asympt-Cn5}.} Following the same argument as in the proof of Theorem \ref{asympt-Cn3-Cn4}, we will only show (\ref{cn5confm}) and (\ref{cn5confl}) can be derived consequently. Using Lemma \ref{EstimatedPath} we have
\begin{align}\label{xbcdist}
\textstyle{\mathbb{P}\Big( \sup\limits_{x\in\mathcal{M}}\Big| \|\wh \X^{\,bc}_x(\wh\theta^{\,bc}_x) -x\| - \frac{|\wh f^{\,bc}(x) - f(x)|}{\|\nabla\wh f^{\,bc}(\wh\X^{\,bc}_x(\wh\theta^{\,bc}_x))\|} \Big| \geq C \tau_n^2 \Big) \leq n^{-L}.}
\end{align}

Similarly, if we consider the trajectories $\wh\X_x$ traveling between $\mathcal{M}^E$ and $\wh{\mathcal{M}}$, then we have
\begin{align*}
\textstyle{\mathbb{P}\Big( \sup\limits_{x\in\mathcal{M}^E}\Big| \|\wh \X_x(\wh\theta_x) -x\| - \frac{|\wh f(x) - \mathbb{E}\wh f(x)|}{\|\nabla\wh f(\wh\X_x(\wh\theta_x))\|} \Big| \geq C \big(\beta_{n,h}^{(0),E}\big)^2\Big) \leq n^{-L}.}
\end{align*}

Let $\mathbb{P}^*$ be the conditional probability measure given $X_1,\cdots,X_n$. Then the bootstrap version of the above result is as follows.
\begin{align}\label{nu0}
\textstyle{\mathbb{E}  \mathbb{P}^*\Big(\sup\limits_{x\in\wh{\mathcal{M}}^{*,E}}\Big| \| \wh\X_x^{*}(\theta_x^{*}) -x\| - \frac{| \widehat f^*(x) - \widehat f^{*,E}(x)|}{\|\nabla \widehat f^{*}(\widehat{\X}_x^{*}(\widehat{\theta}_x^{*}))\|} \Big| \geq C \big(\beta_{n,h}^{(0),E}\big)^2\Big) \leq n^{-L}.}
\end{align}
Notice that $\mathbb{E}\mathbb{P}^*=\mathbb{P}.$ Let $\alpha_n = C_1 \left(\beta_{n,g}^{(0)} + h^2 \right)$ for some $C_1$ fixed and large enough. Let $\mathcal{M}\oplus\alpha_n = \{x\in\mathbb{R}^d: d(x,\mathcal{M})\leq \alpha_n\}$. With $\nu_{n,i}$, $i=1,\cdots, 5$ to be given in (\ref{nu1express}), (\ref{nu2express}), (\ref{nu3express}), (\ref{nu4express}) and (\ref{nu5express}), we will show the following inequalities:
\begin{align}
&\textstyle{\mathbb{P}\Big(\sup\limits_{x\in\mathcal{M}}\Big|  \frac{|\wh f^{\,bc}(x) - f(x)|}{\|\nabla\wh f^{\,bc}(\wh\X^{\,bc}_x(\wh\theta^{\,bc}_x))\|} -  \frac{|\wh f^{\,bc}(x) - f(x)|}{\| \nabla f(x)\|} \Big| \geq C\nu_{n,1}\Big)\leq n^{-L},}\label{nu1}\\
%
&\textstyle{\mathbb{P}\Big( \Big| \sup\limits_{x\in\mathcal{M}\oplus \alpha_n} \frac{|\wh f^{\,bc}(x) - f(x)|}{\| \nabla f(x)\| }  - \sup_{x\in\mathcal{M}} \frac{|\wh f^{\,bc}(x) - f(x)|}{\| \nabla f(x)\| }  \Big| \geq C\nu_{n,2}\Big)\leq n^{-L},}\label{nu2}\\
%
&\textstyle{\mathbb{P}\Big(\Big| \sup\limits_{x\in\mathcal{M}\oplus \alpha_n} \frac{|\wh f^{\,bc}(x) - f(x)|}{\| \nabla f(x)\| }  - \sup\limits_{x\in\mathcal{M}\oplus \alpha_n} \frac{| \widehat f^*(x) - \widehat f^{\,*,E}(x)|}{\| \nabla f(x)\| } \Big| \geq C\nu_{n,3}\Big)\leq n^{-L},}\label{nu3}\\
%
&\textstyle{\mathbb{P}\Big(\Big| \sup\limits_{x\in\mathcal{M}\oplus \alpha_n} \frac{| \widehat f^*(x) - \widehat f^{\,*,E}(x)|}{\| \nabla f(x)\| } - \sup\limits_{x\in\wh{\mathcal{M}}^{*,E}} \frac{| \widehat f^*(x) - \widehat f^{\,*,E}(x)|}{\| \nabla f(x)\| } \Big| \geq C\nu_{n,4}\Big)\leq n^{-L},}\label{nu4}\\
%
%
&\textstyle{\mathbb{P}\Big( \Big| \sup\limits_{x\in\wh{\mathcal{M}}^{*,E}} \frac{| \widehat f^*(x) - \widehat f^{\,*,E}(x)|}{\| \nabla f(x)\| }  - \sup\limits_{x\in\wh{\mathcal{M}}^{*,E}} \frac{| \widehat f^*(x) - \widehat f^{\,*,E}(x)|}{\|\nabla\widehat f^{*}(\widehat{\X}_x^{*}(\widehat{\theta}_x^{*}))\| } \Big| \geq C\nu_{n,5}\Big)\leq n^{-L}.}\label{nu5}
\end{align}


{\em Verification of (\ref{nu1})}. It follows from Lemma~\ref{EstimatedPath} that
\begin{align}\label{xbcdist2}
\sup_{x\in\mathcal{M}} \| \wh\X^{\,bc}_x(\wh\theta^{\,bc}_x) -x \| = O_{a.s.} (\tau_n).
\end{align}
In other words, there exists $C_0>0$ such that $\sup_{x\in\mathcal{M}} \| \wh\X^{\,bc}_x(\wh\theta^{\,bc}_x) -x \|  \leq C_0\tau_n$ for $n$ large enough with probability one. Then we have 
\begin{align*}
&\sup_{x\in\mathcal{M}}\left| \|\nabla\wh f^{\,bc}(\wh\X^{\,bc}_x(\wh\theta^{\,bc}_x))\| - \| \nabla f(x)\| \right|\\
\leq\;  &\sup_{x\in\mathcal{M}}  \|\nabla\wh f^{\,bc}(\wh\X^{\,bc}_x(\wh\theta^{\,bc}_x)) - \nabla f(x)\| \\
\leq\; & \sup_{x\in\mathcal{M},\|y-x\|\leq C_0\tau_n}  \|\nabla\wh f^{\,bc}(y) - \nabla f(x)\|\\
\leq\; &\sup_{y\in\mathcal{M}\oplus (C_0\tau_n)}  \|\nabla\wh f^{\,bc}(y) - \nabla f(y)\| + \sup_{\|y-x\|\leq C_0\tau_n} \|\nabla f(x) - \nabla f(y) \| \\
=\; & O_{a.s.}\left(\beta_{n,h}^{(1),E} + h^2\left(  \beta_{n,l}^{(3),E} + 1 \right) + \tau_n \right),
\end{align*}
where we have used (\ref{Einmahl1}) and (\ref{Einmahl2}). Thus,
\begin{align}\label{nu1express}
\nu_{n,1}=\tau_n \left(\beta_{n,h}^{(1),E} + h^2\left(  \beta_{n,l}^{(3),E} + 1 \right) + \tau_n \right).
\end{align}

{\em Verification of (\ref{nu2})}. We have that

\begin{align*}
&\textstyle{\Big| \sup\limits_{x\in\mathcal{M}\oplus \alpha_n} \frac{|\wh f^{\,bc}(x) - f(x)|}{\| \nabla f(x)\| }  - \sup\limits_{x\in\mathcal{M}} \frac{|\wh f^{\,bc}(x) - f(x)|}{\| \nabla f(x)\| } \Big|}\\
\leq &\textstyle{ \sup\limits_{x\in\mathcal{M}, \|x-y\|\leq\alpha_n} \Big|  \frac{|\wh f^{\,bc}(x) - f(x)|}{\| \nabla f(x)\| } -  \frac{|\wh f^{\,bc}(y) - f(y)|}{\| \nabla f(y)\| } \Big|} \\
\leq& \frac{1}{\epsilon_0^2} \sup_{x\in\mathcal{M}, \|x-y\|\leq\alpha_n} \|\nabla f(x) -\nabla f(y) \| |\wh f^{\,bc}(x) - f(x)|\\
& \hspace*{2cm}+ \frac{1}{\epsilon_0} \sup_{x\in\mathcal{M}, \|x-y\|\leq\alpha_n} |[\wh f^{\,bc}(x) - f(x)] - [\wh f^{\,bc}(y) - f(y)]|.
\end{align*}
Note that
\begin{align*}
\textstyle{\mathbb{P}\Big( \sup\limits_{x\in\mathcal{M}\atop \|x-y\|\leq\alpha_n} \|\nabla f(x) -\nabla f(y) \| |\wh f^{\,bc}(x) - f(x)| \leq C\alpha_n\tau_n\Big)\leq n^{-L},}
\end{align*}
and
\begin{align*}
& \textstyle{\mathbb{P}\Big( \sup\limits_{x\in\mathcal{M} \atop\|x-y\|\leq\alpha_n} |[\wh f^{\,bc}(x) - f(x)] - [\wh f^{\,bc}(y) - f(y)]| 
\geq C\alpha_n\big(\beta_{n,h}^{(1),E} + h^2\big(\beta_{n,l}^{(3),E} + 1\big) \big)\Big)}\\
&\hspace*{3cm}\leq n^{-L},
\end{align*}
where the last equality is obtained following a similar argument on pages 208-209 of Mammen and Polonik (2013). Then we have 
\begin{align}\label{nu2express}
\nu_{n,2} = \alpha_n\big( \beta_{n,h}^{(1),E} + h^2\big(\beta_{n,l}^{(3),E}+1\big) \big).
\end{align}

\noindent
{\em Verification of (\ref{nu3})}. Note that
\begin{align*}
& \textstyle{\Big| \sup\limits_{x\in\mathcal{M}\oplus \alpha_n} \frac{|\wh f^{\,bc}(x) - f(x)|}{\| \nabla f(x)\| }  - \sup\limits_{x\in\mathcal{M}\oplus \alpha_n} \frac{| \widehat f^*(x) - \widehat f^{\,*,E}(x)|}{\| \nabla f(x)\| } \Big|}\\
\leq & \textstyle{\sup\limits_{x\in\mathcal{M}\oplus \alpha_n}  \Big| \frac{|\wh f^{\,bc}(x) - f(x)| - | \widehat f^*(x) - \widehat f^{\,*,E}(x)|}{\| \nabla f(x)\| } \Big| }\\
\leq & \textstyle{\frac{1}{\epsilon_0}\sup\limits_{x\in\mathcal{M}\oplus \alpha_n}  \left| [\wh f(x) - \mathbb{E}\wh f(x)] - [ \widehat f^*(x) - \widehat f^{\,*,E}(x)] \right| + \frac{1}{\epsilon_0}\sup\limits_{x\in\mathcal{M}\oplus \alpha_n} |\wh\beta(x) -\beta(x)|.}
\end{align*}

With (\ref{expcomparison}) and (\ref{2ndderivrate2}), we have 
\begin{align}\label{nu3express}
\textstyle{\nu_{n,3} = \Big(\beta_{n,h}^{(0),E} + \sqrt{\beta_{n,g}^{(0)}}\,\Big)\,\beta_{n,h}^{(0),E} + h^2\big(\beta_{n,l}^{(2)}+h^2\big).}
\end{align}

{\em Verification of (\ref{nu4})}. First observe that
\begin{align*}
\widehat f^{\,*,E}(x) &=  \frac{1}{nh^{d}g^d}\sum_{i=1}^n\int_{\mathbb{R}^d} K\left( \frac{x-y}{h}\right)K\left(\frac{y-X_i}{g}\right)dy\\
&=\frac{1}{ng^d} \sum_{i=1}^n\int_{\mathbb{R}^d} K\left( z\right)K\left(\frac{x-hz-X_i}{g}\right)dz.
%
\end{align*}
We want to show 
\begin{align}\label{bootstrapexpectation}
\mathbb{P}\Big\{ \sup_{x\in\mathbb{R}^d} | \widehat f^{\,*,E}(x)  - f(x)| < C \alpha_n \Big\} = O(n^{-L}).
\end{align}
On the one hand, $\sqrt{n}[\widehat f^{\,*,E}(x) - \mathbb{E} \widehat f^{\,*,E}(x) ]$ can be viewed as an empirical process indexed by
\begin{align*}
\Big\{ \frac{1}{g^d} \int_{\mathbb{R}^d} K( z)K\Big(\frac{x-hz - \cdot}{g}\Big)dz,\;x\in\mathbb{R}^d \Big\}.
\end{align*}
Under our assumptions on the bandwidth and kernel function, following a similar argument on page 207 of Mammen and Polonik (2013), we have
\begin{align}
\mathbb{P}\Big\{ \sup_{x\in\mathbb{R}^d} |\widehat f^{\,*,E}(x) - \mathbb{E} \widehat f^{\,*,E}(x)| > C  \beta_{n,g}^{(0),E}\Big\} = O(n^{-L}).
\end{align}
On the other hand, notice that
\begin{align*}
\mathbb{E} \widehat f^{\,*,E}(x) = \int_{\mathbb{R}^d}\int_{\mathbb{R}^d} K(z) K(u) f(x-hz-gu)dzdu.
\end{align*}
A standard argument leads to 
\begin{align*}
\sup_{x\in\mathbb{R}^d} |\mathbb{E} \widehat f^{\,*,E}(x)  - f(x)| = O(h^2 + g^2).
\end{align*}
Then (\ref{bootstrapexpectation}) is a consequence of the above results. Therefore following the same argument as in Cuevas et al. (2006), we have $d_H(\wh{\mathcal{M}}^{\,*,E}, \mathcal{M}) = O_{a.s.}\left(\alpha_n\right),$ 
which implies that for $C_1$ large enough and when the sample size is large enough, $\wh{\mathcal{M}}^{\,*,E} \subset \mathcal{M} \oplus \alpha_n$ with probability one. 
Therefore
\begin{align*}
& \Big| \sup_{x\in\mathcal{M}\oplus \alpha_n} \frac{| \widehat f^*(x) - \widehat f^{\,*,E}(x)|}{\| \nabla f(x)\| } - \sup_{x\in\wh{\mathcal{M}}^{*,E}} \frac{| \widehat f^*(x) - \widehat f^{\,*,E}(x)|}{\| \nabla f(x)\| } \Big| \\
\leq & \sup_{x\in \wh{\mathcal{M}}^{*,E}, \|x-y\|\leq2\alpha_n} \Big| \frac{| \widehat f^*(x) - \widehat f^{\,*,E}(x)|}{\| \nabla f(x)\| } - \frac{| \widehat f^*(y) - \widehat f^{\,*,E}(y)|}{\| \nabla f(y)\| } \Big|.
\end{align*}
Following the same argument as for (\ref{nu2}), we have
\begin{align*}
\mathbb{E}\mathbb{P}^*\Big(\Big| \sup_{x\in\mathcal{M}\oplus \alpha_n} \frac{| \widehat f^*(x) - \widehat f^{\,*,E}(x)|}{\| \nabla f(x)\| } - \sup_{x\in\wh{\mathcal{M}}^{*,E}} \frac{| \widehat f^*(x) - \widehat f^{\,*,E}(x)|}{\| \nabla f(x)\| } \Big| \geq C\nu_{n,4}\Big)\leq n^{-L}
\end{align*}
with\\[-23pt]
\begin{align}\label{nu4express}
\nu_{n,4} = \alpha_n \beta_{n,h}^{(1),E} .
\end{align}

%


{\em Verification of (\ref{nu5})}. Following an argument similar to (\ref{xbcdist2}), there exists $C_0>0$ such that $\sup_{x\in\wh{\mathcal{M}}^{*,E}} \| \wh\X_x^*(\wh\theta_x^*) -x\|  \leq C_0\beta_{n,h}^{(0),E}$ when the sample size is large enough with probability one. Then we have 
\begin{align*}
&\sup_{x\in\wh{\mathcal{M}}^{*,E}}\left| \|\nabla\wh f^{*}(\wh\X^{*}_x(\wh\theta^{*}_x))\| - \| \nabla f(x)\| \right| \\
&\leq \sup_{x\in\wh{\mathcal{M}}^{*,E}}  \|\nabla\wh f^{*}(\wh\X^{*}_x(\wh\theta^{*}_x)) - \nabla f(x)\| \\
\leq & \sup_{x\in\wh{\mathcal{M}}^{*,E},\|y-x\|\leq C_0\beta_{n,h}^{(0),E}}  \|\nabla\wh f^{*}(y) - \nabla f(x)\|\\
\leq & \sup_{y\in\wh{\mathcal{M}}^{*,E}\oplus \left(C_0\beta_{n,h}^{(0),E}\right)}  \|\nabla\wh f^{*}(y) - \nabla f(y)\| + \sup_{\|y-x\|\leq C_0\beta_n^E} \|\nabla f(x) - \nabla f(y) \| \\
\leq & \sup_{y\in\wh{\mathcal{M}}^{*,E}\oplus \left(C_0\beta_{n,h}^{(0),E}\right)}  \|\nabla\wh f^{*}(y) - \nabla\wh f^{*,E}(y)\|\\
& + \sup_{y\in\wh{\mathcal{M}}^{*,E}\oplus \left(C_0\beta_{n,h}^{(0),E}\right)}  \| \nabla\wh f^{*,E}(y) - \nabla f(y)\| + \sup_{\|y-x\|\leq C_0\beta_n^E} \|\nabla f(x) - \nabla f(y) \|.
%
\end{align*}

Here 
\begin{align*}
&\nabla\widehat f^{*,E}(x) = \frac{1}{nh^{d+1}g^d}\sum_{i=1}^n\int\nabla K\left( \frac{x-y}{h}\right)K\left(\frac{y-X_i}{g}\right)dy\\
=&  \frac{1}{nhg^d}\sum_{i=1}^n\int\nabla K\left( z\right)K\left(\frac{x-hz-X_i}{g}\right)dz\\
=&  \frac{1}{ng^{d+1}}\sum_{i=1}^n\int K\left( z\right)\nabla K\left(\frac{x-hz-X_i}{g}\right)dz.
\end{align*}

Similar to (\ref{bootstrapexpectation}), we have
\begin{align*}
\textstyle{\mathbb{P}\Big\{ \sup\limits_{y\in\wh{\mathcal{M}}^{*,E}\oplus \left(C_0\beta_{n,h}^{(0),E}\right)}  \| \nabla\wh f^{*,E}(y) - \nabla f(y)\| > C \left( \beta_{n,g}^{(1)} + h^2\right) \Big\}= O(n^{-L}).}
\end{align*}

Also, following standard arguments, we obtain
\begin{align*}
\mathbb{E} \mathbb{P}^* \Big\{ \sup_{y\in\wh{\mathcal{M}}^{*,E}\oplus \left(C_0\beta_{n,h}^{(0),E}\right)}  \|\nabla\wh f^{*}(y) - \nabla\wh f^{*,E}(y)\| > C \beta_{n,h}^{(1),E} \Big\}= O(n^{-L}),
\end{align*}
and
\begin{align*}
\sup_{\|y-x\|\leq C_0\beta_{n,h}^{(0),E}} \|\nabla f(x) - \nabla f(y) \| = O\left( \beta_{n,h}^{(0),E} \right).
\end{align*}
Therefore
\begin{align}\label{nu5express}
\nu_{n,5} = \beta_{n,h}^{(0),E} \left( \beta_{n,g}^{(1)} + \beta_{n,h}^{(1)} \right).
\end{align}
After collecting the results (\ref{xbcdist}), (\ref{nu0}), (\ref{nu1}), (\ref{nu2}), (\ref{nu3}), (\ref{nu4}) and (\ref{nu5}), we have
\begin{align}\label{nucombine}
&\mathbb{P}\Big( \Big| \sup_{x\in\mathcal{M}}\|\wh \X^{\,bc}_x(\wh\theta^{\,bc}_x) -x\| - \sup_{x\in\wh{\mathcal{M}}^{*,E}} \| \wh\X_x^{*}(\theta_x^{*}) -x\| \Big| \geq C(\gamma_n^{bc} +\zeta_n/\sqrt{nh^d\log{n}})\Big) \nonumber\\
&\hspace*{3cm}\leq n^{-L},
\end{align}
because $\gamma_n^{bc}+\zeta_n/\sqrt{nh^d\log{n}}$ is the leading term of $\tau_n^2+\big(\beta_{n,h}^{(0),E}\big)^2 + \nu_{n,1} +\nu_{n,2} +\nu_{n,3} +\nu_{n,4} +\nu_{n,5}$.\\[5pt]
Similar to (\ref{neumannres4}), we obtain
\begin{align}\label{neumannres5}
\textstyle{\mathbb{P} \Big( \sup\limits_{x\in\mathcal{M}} \frac{|\wh f^{\,bc}(x) - f(x) |}{\|\nabla f(x)\|} \in [c, d]\Big) = O\big( \Psi_n(d-c)+ h^2\big(\beta_{n,l}^{(2)}+h^2\big)\sqrt{nh^d\log{n}}\big).}
\end{align}
Note that the proof of the above result need to adapt the proof of Proposition 3.1 in Neumann (1998). Specifically, using Neumann's notation, corresponding to page 2045 in Neumann (1998), we have that 
\begin{align*}
\sup_{k\in\mathscr{K}_l}\Big\{ \sup_{x\in I_k} \Big(\frac{T_{k1}}{\|\nabla f(x)\|} + \frac{T_{k2}(x)}{\|\nabla f(x)\|}\Big) \Big\} \in [c,d]
\end{align*}
implies 
\begin{align*}
\sup_{k\in\mathscr{K}_l}\Big\{  \frac{T_{k1}}{\|\nabla f(x_k)\|} + \sup_{x\in I_k} \frac{T_{k2}(x)}{\|\nabla f(x)\|} \Big\} \in [c-\kappa_n,d+\kappa_n],
\end{align*}
where $x_k$ is a fixed point on $I_k$, say $(2(k_1-1)h,\cdots,2(k_d-1)h)$, and $\kappa_n = Ch \beta_{n,h}^{(0),E}$. This is because $I_k$ is a cube with size length of $h$, $\|\nabla f\|$ is differentiable and bounded away from zero in a neighborhood of the level set. Then the rest of the proof follows the proof of Proposition 3.1 in Neumann (1998).\\

With (\ref{xbcdist}), (\ref{nu1}) and (\ref{neumannres5}), we have with $\lambda_n=C(\tau_n^2 + \nu_{n,1})$,
\begin{align}\label{neumannres6}
&\mathbb{P} \Big( \sup_{x\in\mathcal{M}}\|\wh \X^{\,bc}_x(\wh\theta^{\,bc}_x) -x\| \in [c, d]\Big) \nonumber\\
\leq\; & \mathbb{P} \Big( \sup_{x\in\mathcal{M}} \frac{|\wh f^{\,bc}(x) - f(x) |}{\|\nabla f(x)\|}  \in [c-\lambda_n, d+\lambda_n]\Big) + O(n^{-L}) \nonumber\\
=\;& O\Big( \Psi_n(d-c)+ h^2\big(\beta_{n,l}^{(2)}+h^2\big)\sqrt{nh^d\log{n}} + \lambda_n \sqrt{nh^d\log{n}}\Big).
\end{align}
An application of Lemma~\ref{MP} with $Z_n = \sup_{x\in\mathcal{M}}\|\wh \X^{\,bc}_x(\wh\theta^{\,bc}_x) -x\| $ and $Z_n^* = \sup_{x\in\wh{\mathcal{M}}^{*,E}} \| \wh\X_x^{*}(\theta_x^{*}) -x\| $, by using (\ref{nucombine}) and (\ref{neumannres6}), concludes the proof.\hfill$\square$

\section{Appendix}
The following result is essentially taken from Mammen and Polonik (2013). We state it here for easy reference.\\[5pt]
Using the notation introduced above, let $Z_n$ be a statistic, and let $Z_n^*$ be a bootstrap version of this statistic.   
For $0 < \alpha < 1$, define
\begin{align*}
%
&\wh c^*_n(1-\alpha) = \sup\big\{t: P^*\big(Z_n^* \le t) \le 1-\alpha   \big\}
\end{align*}
and let $c_n(1-\alpha)$ be defined similarly with $Z_n^*$ replaced by $Z_n$ (and $P^*$ replaced by $P$). 

\begin{lemma}\label{MP}
Suppose that there exist sequences $\{\gamma_n\}, \{\delta_n\}$ and $\{\tau_n\}$ such that
\begin{align}
 &\mathbb{P}\big(|Z_n - Z_n^*| > \gamma_n\big) \le \delta_n\quad \text{and} \label{cond1}\\
  \sup_{t \in \R}\,&\mathbb{P}\big(Z_n \in [t,t+\gamma_n)\big) \le \tau_n.\label{cond2}
\end{align}
Then we have, for $0 < \alpha < 1,$

\begin{align}\label{rate}
\big|\mathbb{P}\big( Z_n \le \wh c^*_n(1-\alpha)\big) - (1 - \alpha)\big| \le  7\tau_n + 5\sqrt{\delta_n}.
\end{align}
\end{lemma}
{\bf Proof.} Lemma 2.4 of Mammen and Polonik (2013) says that under the stated conditions, we have
\begin{align*}
    \big| \mathbb{P}\big(Z_n \le \wh c^*_n(1-\alpha)\big) - \mathbb{P}\big(Z_n \le c_n(1-\alpha)\big) \big| \le 6 \tau_n + 5 \sqrt{\delta_n\,}.
\end{align*}
  It remains to observe that, by using assumption (\ref{cond2}),
\begin{align*}
    \big|\mathbb{P}\big(Z_n \le c_n(1-\alpha)\big) - (1-\alpha)\big| 
    & \le \sup_{t \in \R}\mathbb{P}\big(Z_n \in [t-\gamma_n,t+\gamma_n)\big) \le \tau_n.
\end{align*}
Note that in Mammen and Polonik (2013), the quantities $Z_n$ and $Z_n^*$ denote particular statistics, their Lemmas 2.2 and 2.4 in fact hold for any statistics satisfying (\ref{cond1}) and (\ref{cond2}). Indeed, an inspection of their proofs shows that they do not use any other properties of the statistics. 
\hfill $\square$\\

\section*{Acknowledgement}
The authors would like to thank the Associate Editor and the referee for careful reading of the manuscript and for insightful comments that lead to significant improvements. The research of Wolfgang Polonik was partially supported by NSF grant DMS 1107206. The research of Wanli Qiao was partially supported by NSF grant DMS 1821154. The simulations in this work were run on ARGO, a research computing cluster provided by the Office of Research Computing at George Mason University, VA.

\section*{References}

\begin{description}
\item Ambrosio, L., Colesanti, A. and Villa, E. (2008): Outer Minkowski content for some classes of closed sets. {\em Math. Ann.} {\bf 342}, 727-748.
\item Arias-Castro, E., Mason, D., and Pelletier, B. (2016): On the estimation of the gradient lines of a density and the consistency of the mean-shift algorithm. {Journal of Machine Learning Research} {\bf 17} 1-28.
\item Audibert, J-Y. and Tsybakov, A. (2007): Fast Learning rates for plug-in classifier. {\em Ann. Statist} {\bf 35},608-633. 
\item Biau, G., Cadre, B., and Pelletier, B. (2007): A graph-based estimator of the number of clusters. {\em ESAIM Probab. Stat.} {\bf 11} 272-280.
\item Bickel, P. and Rosenblatt, M. (1973): On some global measures of the deviations of density function estimates. {\em Ann. Statist.}, {\bf 1}(6), 1071-1095.
\item Bobrowski, O., Mukherjee, S. and Taylor, J.E. (2017): Topological consistency via kernel estimation. {\em Bernoulli}, {\bf 23}, 288 - 328.
\item Bredon, G.E. (1993): Topology and Geometry. Volume 139 of {\em Graduate Texts in Mathematics}. Springer-Verlag, New York.
\item Breuning, M.M., Kriegel, H.P., Ng R.T., and Sander, J. (2000): Lof: identifying density-based local outlier. {\em ACM sigmod record}, {\bf 29}, 93-104.
\item Broida, J.G. and Willamson, S.G. (1989): {\em A Comprehensive Introduction to Linear Algebra.} Addison-Wesley
\item Cadre, B. (2006): Kernel estimation of density level sets. {\em J. Multivariate Anal.} {\bf 97} 999-1023.
\item Calonico, S., Cattaneo, M.D. and Farrell, M.H. (2018a): On the effect of bias estimation on coverage accuracy in nonparametric inference. {\em Journal of the American Statistical Association}, DOI: 10.1080/01621459.2017.1285776.
\item Calonico, S., Cattaneo, M.D. and Farrell, M.H. (2018b): Coverage error optimal confidence intervals. {\em arXiv:1808.01398}
\item Cavalier, L (1997): Nonparametric estimation of regression level sets. {\em Statistics.} {\bf 29}, 131-160.
\item Chazal, F., Lieutier, A. and Rossignac, J. (2007): Normal-map between normal-compatible manifolds. {\em International Journal of
Computational Geometry and and Applications}, {\bf 17}, 403-421.
\item Chen, Y, Genovese, C.R., Wasserman, L (2017): Density level set: asymptotics, inference, and visualization. {\em J. Amer. Statist. Assoc.}, {\bf 112} 1684-1696.
\item Chen, Y. (2017): Nonparametric Inference via Bootstrapping the Debiased Estimator. {\em arXiv: 1702.07027}
\item Chernozhukov, V., Chetverikov, D. and Kato, K. (2014): Gaussian approximation of suprema of empirical processes. {\em Ann. Statist} {\bf 42}, 1564-1597.
\item Cuevas, A (2009): Set estimation: Another bridge between statistics and geometry. {\em Bolet\'{i}n de Estad\'{i}stica e Investigaci\'{o}n Operativa}.
\item Cuevas, A., Febrero, M. and Fraiman, R. (2000): Estimating the number of clusters. {\em Canad. J. Statist.} {\bf 28}, 367-382.
\item Cuevas, A., Fraiman, R., and Pateiro-L\'{o}pez, B. (2012): On statistical properties of sets fulfilling rolling-type conditions. {\em Advances in Applied Probability} {\bf 44} 311-329.
\item Cuevas, A., Gonz\'{a}lez-Manteiga, W., and  Rodr\'{i}guez-Casal, A. (2006): Plug-in estimation of general level sets. {\em Australian \& New Zealand Journal of Statistics} {\bf 48} 7-19.
\item Cuevas, A. and Rodr\'{i}guez-Casal, A. (2004): On boundary estimation. {\em Advances in Applied Probability}, 340-354.
\item Einmahl, U., and Mason, D.M. (2005): Uniform in bandwidth consistency of kernel-type function estimators. {\em Ann. Statist.}, {\bf 33}, 1380-1403.
\item Fasy, B.T., Lecci, F., Rinaldo, A., Wasserman, L., Balakrishnan, S., and Singh, L. (2014): Confidence sets for persistence diagrams. {\em Ann. Statist.}, {\bf 42}, 2301-2339.
\item Federer, H. (1959): Curvature measures. {\em Transactions of the American Mathematical Society}, {\bf 93}, 418-491.
%
%
\item Hall, P. (1979): The rate of convergence of normal extremes. {J. Appl. Probab.}, {\bf 16}, 433-439.
\item Hall, P. (1992): {\em The Bootstrap and Edgeworth Expansion}. Springer-Verlag, New York.
\item Hall, P. (1993): On Edgeworth expansion and bootstrap confidence bands in nonparametric curve estimation. {\em Journal of the Royal Statistical Society, Series B}. {\bf 55}, 291-304.
\item Hall, P. and Jing, B.-Y. (1995): Uniform Coverage Error Bounds for Confidence Intervals and Berry-Esseen Theorems for Edgeworth Expansion. {\em Annals of Statistics}, {\bf 23}, 363-375.
\item Hall, P. and Kang, K-H. (2005): Bandwidth choice for nonparametric classification. {\em Ann. Statist}. {\bf 33}, 284-306.
\item Hartigan, J.A. (1987): Estimation of a convex density contour in two dimensions. {\em J. Amer. Statist. Assoc.}, {\bf 82}, 267 - 270.
\item Hodge, V.J., and Austin, J. (2004): A survey of outlier detection methodologies. {\em Artificial Intelligence Review}, {\bf 22}(2), 85-126.
\item Jang, W. (2006): Nonparametric density estimation and clustering in astronomical sky survey. {\em Comp. Statist. \& Data Anal.} {\bf 50}, 760 - 774.
\item Jankowski, H and Stanberry, L. (2014): Visualizing variability: Confidence regions in level set estimation. Proceedings of the 16th International Conference on Geometry and Graphics, 1328-1339.
%
%
\item  Mammen, E. and Polonik, W. (2013): Confidence sets for level sets. {\em Journal of Multivariate Analysis}, {\bf 122}(C), 202-214.
\item Mammen, E. and Tsybakov, A.B. (1999): Smooth discrimination analysis. {\em Ann. Statist.} {\bf 27}, 1808-1829.
\item Mason, D. and Polonik, W. (2009): Asymptotic normality of plug-in level set estimates {\em Annals of Applied Probability}, {\bf 19}(3), 1108-1142.
\item Neumann, M.H. (1998): Strong approximation of density estimators from weakly dependent observations by density estimators from independent observations. {\em Ann. Statist.} {\bf 26}, 2014-2048.
\item  Piterbarg, V.I. (1996): \emph{Asymptotic Methods in the Theory of Gaussian Processes and Fields}, Translations of Mathematical Monographs, Vol. 148, American Mathematical Society, Providence, RI.
\item Polonik, W. (1995): Measuring mass concentrations and estimating density contour clusters - an excess mass approach. {\em Ann. Statist.} {\bf 23}, 855-881.
\item Qiao, W. (2018a): Asymptotics and optimal bandwidth selection for nonparametric estimation of density level sets. {\em arXiv: 1707.09697}.
\item Qiao, W. (2018b): Nonparametric estimation of surface integrals on density level sets. {\em arXiv: 1804.03601}.
\item Qiao, W. and Polonik, W. (2018): Extrema of rescaled locally stationary Gaussian fields on manifolds. {\em Bernoulli} {\bf 24}, 1834-1859.
\item Qiao, W. and Polonik, W. (2016): Theoretical analysis of nonparametric filament estimation. {\em Ann. Statist.} {\bf 44}, 1269-1297.
\item Rinaldo, A., Singh, A., Nugent, R. and Wasserman, L. (2010): Stability of density-based clustering. {\em arXiv: 1011.2771v1}.
\item Rosenblatt, M. (1976): On the maximal deviation of $k$-dimensional density estimates. {\em Ann. Probab.}, {\bf 4}(6), 1009-1015.
\item Samworth, R.J. and Wand, M.P. (2010): Asymptotics and optimal bandwidth selection for highest density region estimation. {\em Ann. Statist.} {\bf 38} 1767-1792.
\item Sommerfeld, M., Sain, S., and Schwartzman, A. (2015): Confidence regions for excursion sets in asymptotically Gaussian random fields, with an application to climate. {\em arXiv: 1501.07000}.
\item Steinwart, I., Hush, D. and Scovel, C. (2005): A classification framework for anomaly detection. {\em J. Machine Learning Reserach} {\bf 6}, 211-232.
\item Tsybakov, A.B. (1997): Nonparametric estimation of density level sets. {\em Ann. Statist.} {\bf 25}, 948-969.
\item Walther, G. (1997): Ganulometric smoothing. {\em Ann. Statist.} {\bf 25}, 2273- 2299. 
\item Willett, R.M. and Nowak, R.D. (2005): Level set estimation in medical imaging, {\em Proceedings of the IEEE Statistical Signal Processing, Vol. 5}, 1089-1092.
\end{description}

\end{document}